\documentclass[11pt]{article} \usepackage{fullpage}

\usepackage[sc]{mathpazo}
\usepackage[T1]{fontenc}
\usepackage{microtype}
\usepackage{abstract}

\usepackage{amsfonts}

\usepackage{graphicx} %
\usepackage{subfigure} \usepackage{color}

\usepackage{enumitem}
\usepackage{wrapfig}
\usepackage{booktabs}
\usepackage{natbib}

\usepackage{algorithm} \usepackage{algorithmic}

\usepackage{epsfig} \usepackage{amsthm} \usepackage{amsmath}
\usepackage{amssymb}
\usepackage{color}
\usepackage{epstopdf} \usepackage{hyperref}

\usepackage{breqn}

\newcounter{ass_counter} \newcounter{thm_counter}

\newtheorem{theorem}[thm_counter]{Theorem}%
\newtheorem{corollary}[thm_counter]{Corollary}
\newtheorem{assumption}[ass_counter]{Assumption}

\def\E{\mathbb{E}}

 \newcommand{\AsySG}{\mbox{\sc
    AsySG}} \newcommand{\AsySGC}{\mbox{\sc AsySG-con}}
\newcommand{\AsySGI}{\mbox{\sc AsySG-incon}}
\newcommand{\hogwild}{\mbox{\sc Hogwild!}}

\newcommand{\Lmax}{L_{\mbox{\rm\scriptsize max}}}
\newcommand{\beq}{\begin{equation}} \newcommand{\eeq}{\end{equation}}

\allowdisplaybreaks[4]
\title{Asynchronous Parallel Stochastic Gradient for Nonconvex
Optimization\footnote{Published at NIPS 2015.}}

\author{Xiangru Lian, Yijun Huang, Yuncheng Li, and Ji Liu\\
\{\texttt{lianxiangru, huangyj0, raingomm, ji.liu.uwisc\}@gmail.com}\\
Department of Computer Science, University of Rochester}
\date{}

\begin{document}
\maketitle %

\begin{abstract}
  Asynchronous parallel implementations of stochastic gradient (SG)
  have been broadly used in solving deep neural network and received
  many successes in practice recently. However, existing theories
  cannot explain their convergence and speedup properties, mainly due
  to the nonconvexity of most deep learning formulations and the
  asynchronous parallel mechanism. To fill the gaps in theory and
  provide theoretical supports, this paper studies two asynchronous
  parallel implementations of SG: one is  over a
  computer network and the other is on  a shared
  memory system. We establish an ergodic convergence rate
  $O(1/\sqrt{K})$ for both algorithms and prove that the linear
  speedup is achievable if the number of workers is bounded by
  $\sqrt{K}$ ($K$ is the total number of iterations). Our results
  generalize and improve existing analysis for convex minimization.

\end{abstract}

\section{Introduction}
The asynchronous parallel optimization recently received many successes and
broad attention in machine learning and optimization \citep{recht2011hogwild,
  li2013parameter, li2014communication, yun2013nomad, fercoq2013accelerated,
  zhang2014asynchronous, marecek2014distributed, tappenden2015complexity,
  hong2014distributed}. It is mainly due to that the asynchronous parallelism
largely reduces the system overhead comparing to the synchronous parallelism.
The key idea of the asynchronous parallelism is to allow all workers work
independently and have no need of synchronization or coordination. The
asynchronous parallelism has been successfully applied to speedup many
state-of-the-art optimization algorithms including stochastic gradient
\citep{recht2011hogwild, agarwal2011distributed, ZhangCL14a,
  2015arXiv150504824F, paine2013gpu, mania2015perturbed}, stochastic coordinate
descent \citep{avron2014revisiting, liu2013asynchronous,
  sridhar2013approximate}, dual stochastic coordinate ascent
\citep{asySDCicml2015}, and randomized Kaczmarz algorithm
\citep{liu2014asynchronousRK}.

In this paper, we are particularly interested in the asynchronous
parallel stochastic gradient algorithm ($\AsySG$) for \emph{nonconvex}
optimization mainly due to its recent successes and popularity in deep
neural network \citep{bengio2003neural, dean2012large, paine2013gpu,
  ZhangCL14a, li2014scaling} and matrix completion
\citep{recht2011hogwild, petroni2014gasgd, yun2013nomad}. While some
research efforts have been made to study the convergence and speedup
properties of $\AsySG$ for \emph{convex} optimization, people still
know very little about its properties in \emph{nonconvex}
optimization. Existing theories cannot explain its convergence and
excellent speedup property in practice, mainly due to the nonconvexity
of most deep learning formulations and the asynchronous parallel
mechanism. People even have no idea if its convergence is certified
for nonconvex optimization, although it has been used widely in
solving deep neural network and implemented on different platforms
such as computer network and shared memory (for example, multicore and
multiGPU) system.

To fill these gaps in theory, this paper tries to make the first attempt
to study $\AsySG$ for the following nonconvex optimization problem
\begin{equation}
  \min_{x\in \mathbb{R}^n}\quad  f(x) := \mathbb{E}_\xi[F(x; \xi )]
  \label{eq:main}
\end{equation}
where $\xi \in \Xi$ is a random variable and $f(x)$ is a smooth (but
not necessarily convex) function. The most common specification is
that $\Xi$ is an index set of all training samples
$\Xi = \{1, 2, \cdots, N\}$ and $F(x; \xi)$ is the loss function with
respect to the training sample indexed by $\xi$.

We consider two popular asynchronous parallel implementations of SG:
one is for the computer network originally proposed in
\citep{agarwal2011distributed} and the other one is for the shared
memory (including multicore/multiGPU) system originally proposed in
\citep{recht2011hogwild}. Note that due to the architecture diversity,
it leads to two different algorithms. The key difference lies on that
the computer network can naturally (also efficiently) ensure the
atomicity of reading and writing the \emph{whole} vector of $x$,
while the shared memory system is unable to do that efficiently and
usually only ensures efficiency for atomic reading and writing on a
\emph{single} coordinate of parameter $x$.
The implementation on computer cluster is described by the
``consistent asynchronous parallel SG'' algorithm ($\AsySGC$),
because the value of  parameter $x$
used for stochastic gradient evaluation is {\bf con}sistent -- an
existing value of parameter $x$
at some time point. Contrarily, we use the ``inconsistent asynchronous
parallel SG'' algorithm ($\AsySGI$)
to describe the implementation on the shared memory platform, because
the value of parameter $x$
used is {\bf incon}consistent, that is, it might not be the real state
of $x$ at any time point.

This paper studies the theoretical convergence and speedup properties
for both algorithms. We establish an asymptotic convergence rate
of $O(1/\sqrt{KM})$
for $\AsySGC$
where $K$
is the total iteration number and $M$
is the size of minibatch. The linear speedup\footnote{The speedup for
  $T$
  workers is defined as the ratio between the total work load using
  one worker and the average work load using $T$
  workers to obtain a solution at the same precision. ``The linear
  speedup is achieved'' means that the speedup with $T$
  workers greater than $cT$
  for any values of $T$
  ($c\in (0, 1]$
  is a constant independent to $T$).}
is proved to be achievable while the number of
workers is bounded by $O(\sqrt{K})$.
For $\AsySGI$,
we establish  an asymptotic convergence and
speedup properties similar to $\AsySGC$.
The intuition of the linear speedup of
asynchronous parallelism for SG can be explained in the following: Recall that
the serial SG essentially uses the ``stochastic'' gradient to
surrogate the accurate gradient. $\AsySG$
brings additional deviation from the accurate gradient due to
using ``stale'' (or delayed) information. If the additional deviation
is relatively minor to the deviation caused by the
``stochastic'' in SG, the total iteration complexity (or convergence
rate) of $\AsySG$ would be comparable to the serial SG, which implies
a nearly linear speedup. This is the key reason
why $\AsySG$ works.

The main contributions of this paper are highlighted as follows:
\begin{itemize}
\item Our result for $\AsySGC$
  generalizes and improves earlier analysis of
  $\AsySGC$
  for convex optimization in
  \citep{agarwal2011distributed}. Particularly, we improve the upper
  bound of the maximal number of workers to ensure the linear speedup
  from $O(K^{1/4}M^{-3/4})$
  to $O(K^{1/2}M^{-1/2})$ by a factor $K^{1/4}M^{1/4}$;
\item The proposed $\AsySGI$
  algorithm provides a more accurate description than $\hogwild$
  \citep{recht2011hogwild} for the lock-free
  implementation of $\AsySG$
  on the shared memory system. Although our result does not strictly
  dominate the result for $\hogwild$
  due to different problem settings, our result can be applied
  to more scenarios (e.g., nonconvex
  optimization);
\item Our analysis provides theoretical
  (convergence and speedup) guarantees for many recent successes of
  $\AsySG$
  in deep learning. To the best of our knowledge, this is the first
  work that offers such theoretical support.
\end{itemize}

\subsection*{Notation}
\begin{itemize}
\item $x^*$ denotes the global optimal solution to \eqref{eq:main}.
\item $\|x\|_0$
  denotes the $\ell_0$
  norm of vector $x$, that is, the number of nonzeros in $x$;
\item $e_i \in \mathbb{R}^n$
  denotes the $i$th natural unit basis vector;
\item We use $\mathbb{E}_{\xi_{k, *}}(\cdot)$
  to denote the expectation with respect to a set of variables
  $\{\xi_{k, 1}, \cdots, \xi_{k,M}\}$;
\item $\E(\cdot)$
  means taking the expectation in terms of all random
  variables. %
\item $G(x; \xi)$ is used to denote $\nabla F(x; \xi)$ for short;
\item We use $\nabla_i f(x)$
  and $(G(x; \xi))_i$
  to denote the $i$th
  element of $\nabla f(x)$ and $G(x; \xi)$ respectively.
\end{itemize}

{\noindent \bf Assumption} Throughout this paper, we make the
following assumption for the objective function. All
of them are quite common in
the analysis of stochastic gradient algorithms.
\begin{assumption} \label{ass:0} We assume that the following holds:
  \begin{itemize}
  \item {\bf (Unbiased Gradient):} The stochastic gradient $G(x; \xi)$
    is unbiased, that is to say, \begin{equation} \nabla f(x) =
      \mathbb{E}_\xi[G(x; \xi)]
      \label{ass:1.2}
    \end{equation}
  \item {\bf (Bounded Variance):} The variance of stochastic gradient
    is bounded:
    \begin{equation}
      \mathbb{E}_\xi(\|G(x; \xi) -\nabla f(x)\|^2) \le \sigma^2, \quad \forall x.
      \label{eq:boundedvar}
    \end{equation}
  \item {\bf (Lipschitzian Gradient):} The gradient function
    $\nabla f(\cdot)$ is Lipschitzian, that is to say,
    \begin{equation}
      \|\nabla f(x) - \nabla f(y)\| \le L\|x-y\| \quad \forall x, \forall y.
      \label{eq:Lip1}
    \end{equation}
    Under the Lipschitzian gradient assumption, we can define two more
    constants $L_s$ and $\Lmax$. Let $s$ be any positive
    integer. Define $L_s$ to be the minimal constant satisfying the
    following inequality:
    \begin{equation}
      \left\|\nabla f\left(x\right)- \nabla f\left(x+\sum_{i\in S} \alpha_i e_ i\right) \right\| \leq L_s \left\|\sum\limits_{i\in S} \alpha_i e_i\right\| \condition{$\forall S \subset \{1, 2, ..., n\}$ and $|S| \leq s$}
      \label{ass:algo2-lip}
    \end{equation}
    Define $\Lmax$ as the minimum constant that satisfies:
    \begin{equation}
      |\nabla_i f(x) - \nabla_i f(x+ \alpha e_i)| \le \Lmax|\alpha| \condition{$\forall i \in \{1, 2, ..., n\}$}.
    \end{equation}
    It can be seen that $\Lmax \leq L_s \leq L$.
  \end{itemize}
\end{assumption}

\section{Related Work}
This section reviews stochastic gradient algorithms, synchronous
parallel stochastic gradient algorithms,
asynchronous parallel  gradient algorithms, and
asynchronous parallel stochastic gradient algorithms.

\emph{Stochastic gradient} is a very powerful
approach for large scale optimization. The well known convergence rate
of stochastic gradient is $O(1/\sqrt{K})$
for convex problems and $O(1/{K})$
for strongly convex problems, for example, see
\citep{nemirovski2009robust, moulines2011non}. A parallel algorithm
 of stochastic gradient is the dual averaging
approach, which achieves the same convergence rate
\citep{xiao2009dual}. The stochastic gradient is also closely related
to online learning algorithms. Please refer to the online learning
literature, for example, \citep{crammer2006online,
  shalev2011online, yang2014regret}. Most studies for stochastic
gradient focus on convex optimization. For nonconvex optimization
using SG, \citet{ghadimi2013stochastic}  proved
an ergodic convergence rate of
$O({1/ \sqrt{K}})$
for stochastic gradient under nonconvex optimization, which is
consistent with the rate of stochastic gradient for convex
problems.

As for \emph{synchronous parallel stochastic gradient methods},
\citet{dekel2012optimal} proposed a mini-batch stochastic gradient
algorithm -- multiple workers compute the stochastic gradient on $M$
data in parallel and need a synchronization step before modifying
parameter $x$
in every iteration. The convergence rate is proven to be
$O(1 / \sqrt{KM})$
for convex optimization. In their follow up work
\citep{dekel2010robust}, they showed that a variant of their approach
is asymptotically robust to asynchrony, as long as most processors
remain synchronized for most of time. (Precisely, this variant can be
considered as a partially asynchronous parallel stochastic gradient
algorithm.) \citep{zinkevich2010parallelized} studied a method where
each node in the network runs the vanilla
stochastic gradient method, using random subsets of the overall data
set, and their solutions are only averaged in the
final
step. %

The \emph{asynchronous parallel algorithms}
received broad attention in optimization recently, although pioneer
studies started from 1980s \citep{bertsekas1989parallel}. Due to the
rapid development of hardware resources, the asynchronous parallelism
recently received many successes when applied to
parallel stochastic gradient
\citep{recht2011hogwild, agarwal2011distributed, ZhangCL14a,
  2015arXiv150504824F, paine2013gpu}, stochastic coordinate descent
\citep{avron2014revisiting, liu2013asynchronous}, dual stochastic
coordinate ascent \citep{asySDCicml2015}, randomized Kaczmarz
algorithm \citep{liu2014asynchronousRK}, and ADMM
\citep{zhang2014asynchronous}. \citet{liu2013asynchronous} and
\citet{liu2014asynchronous} studied the asynchronous parallel
stochastic coordinate descent algorithm with consistent read and
inconsistent read respectively and prove the linear speedup is
achievable if $T\leq O(n^{1/2})$
for smooth convex functions and $T\leq O(n^{1/4})$
for functions with ``smooth convex loss +
nonsmooth convex separable
regularization''. \citet{avron2014revisiting} studied this
asynchronous parallel stochastic coordinate descent algorithm
in solving $Ax = b$
where $A$
is a symmetric positive definite matrix, and showed
that the linear speedup is achievable if $T \leq O(n)$
for consistent read and $T \leq O(n^{1/2})$ for inconsistent read.  \citet{asySDCicml2015} studied a semi-asynchronous parallel version
of Stochastic Dual Coordinate Ascent algorithm which periodically
enforces primal-dual synchronization in a separate thread.

We review the \emph{asynchronous parallel stochastic gradient
  algorithms} in the last. \citet{agarwal2011distributed} analyzed the
$\AsySGC$ algorithm (on computer cluster) for \emph{convex smooth}
optimization and proved
a convergence rate of
$O\left({1\over \sqrt{MK}} + {MT^2 \over K}\right)$
which implies that linear speedup is achieved when
$T$
is bounded by $O(K^{1/4}/M^{3/4})$.
In comparison, our analysis for the more general nonconvex smooth
optimization improves the upper bound by a factor $K^{1/4}M^{1/4}$.
A very recent work \citep{2015arXiv150504824F} extended the analysis
in \citet{agarwal2011distributed} to minimize functions in the form
``smooth convex loss + nonsmooth convex regularization'' and obtained
similar results.
\citet{recht2011hogwild} proposed a lock free asynchronous parallel
implementation of SG on the shared memory system and described this
implementation as $\hogwild$ algorithm.
They proved a sublinear convergence rate $O(1/K)$
for \emph{strongly} convex smooth
objectives. Another recent work
\cite{mania2015perturbed} analyzed asynchronous stochastic
optimization algorithms for convex functions by viewing it as a
serial algorithm with the input perturbed by bounded noise and
proved the convergences rates no worse than using traditional point
of view for several algorithms.

\section{Asynchronous parallel stochastic gradient for computer
  network} \label{sec:AsySGC}

This section considers the asynchronous parallel implementation of SG
on computer network proposed by \citet{agarwal2011distributed}. It has
been successfully applied to the distributed neural network
\citep{dean2012large} and the parameter server \citep{li2014scaling}
to solve deep neural network.

\subsection{Algorithm Description: $\AsySGC$}

\begin{algorithm}[H] %
  \caption{$\AsySGC$} %
  \label{alg:conread} %
  \begin{algorithmic} [1] %
    \REQUIRE $x_0$, $K$, $\{\gamma_k\}_{k=0, \cdots, K-1}$ \ENSURE
    $x_K$ \FOR{$k=0,\cdots, K-1$} \STATE Randomly select $M$ training
    samples indexed by $\xi_{k,1}, \xi_{k,2}, ... \xi_{k,M}$; \STATE
    $x_{k+1} = x_{k} - \gamma_k \sum_{m=1}^{M}G(x_{k-\tau_{k,m}},
    \xi_{k,m})$;
    \ENDFOR
  \end{algorithmic}
\end{algorithm}

The ``star'' in the star-shaped network is a master
machine\footnote{There could be more than one
  machines in some networks, but all of them serves the same
  purpose and can be treated as a single machine.}
which maintains the parameter $x$.
Other machines in the computer network serve as workers which only
communicate with the master. All workers exchange information with the
master independently and simultaneously, basically repeating the
following steps:
\begin{itemize}
\item {\bf (Select):} randomly select a subset of training samples
  $S \in \Xi$;
\item {\bf (Pull):} pull parameter $x$ from the master;
\item {\bf (Compute):} compute the stochastic gradient
  $g\leftarrow \sum_{\xi \in S} G(x; \xi)$;
\item {\bf (Push):} push $g$ to the master.
\end{itemize}
The master basically repeats the following steps:
\begin{itemize}
\item {\bf (Aggregate):} aggregate a certain amount of stochastic
  gradients ``$g$'' from workers;
\item {\bf (Sum):} summarize all ``$g$''s into a vector $\Delta$;
\item {\bf (Update):} update parameter $x$ by
  $x \leftarrow x - \gamma \Delta$.
\end{itemize}
While the master is aggregating stochastic gradients from workers, it
does not care about the sources of the collected stochastic
gradients. As long as the total amount achieves the predefined
quantity, the master will compute $\Delta$ and perform the update on
$x$. The ``update'' step is performed as an atomic operation --
workers cannot read the value of $x$ during this step, which can be
efficiently implemented in the network (especially in the parameter
server \citep{li2014scaling}). The key difference between this
asynchronous parallel implementation of SG and the serial (or
synchronous parallel) SG algorithm lies on that in the ``update''
step, some stochastic gradients ``$g$'' in ``$\Delta$'' might be
computed from some early value of $x$ instead of the current one,
while in the serial SG, all $g$'s are guaranteed to use the current
value of $x$.

The asynchronous parallel implementation substantially reduces the
system overhead and overcomes the possible large network delay, but
the cost is to use the old value of ``$x$''
in the stochastic gradient evaluation. We will show
in Section~\ref{sec:con:thm} that the negative
affect of this cost will vanish asymptotically.

To mathematically characterize this asynchronous parallel
implementation, we monitor parameter $x$
in the master. We use the subscript $k$
to indicate the $k$th
iteration on the master. For example, $x_k$
denotes the value of parameter $x$
after $k$
updates, so on and so forth. We introduce a variable $\tau_{k, m}$
to denote how many delays for $x$
used in evaluating the $m$th
stochastic gradient at the $k$th iteration.
This asynchronous parallel implementation of SG on the ``star-shaped''
network is summarized by the $\AsySGC$
algorithm, see Algorithm~\ref{alg:conread}. The suffix
``${\mbox{\sc con}}$''
is short for ``consistent read''. ``Consistent read'' means that the
value of $x$
used to compute the stochastic gradient is a real state of $x$
no matter at which time point. ``Consistent read'' is ensured by the
atomicity of the ``update'' step. When the atomicity fails, it leads
to ``inconsistent read'' which will be discussed in
Section~\ref{sec:asysgi}. It is worth noting that on some ``non-star''
structures the asynchronous implementation can also be described as
$\AsySGC$
in Algorithm~\ref{alg:conread}, for example, the cyclic delayed
architecture and the locally averaged delayed architecture
\citep[Figure 2]{agarwal2011distributed} .

\subsection{Analysis for $\AsySGC$}\label{sec:con:thm}

To analyze Algorithm~\ref{alg:conread}, besides
Assumption~\ref{ass:0} we make the following additional assumptions.
\begin{assumption} \label{ass:1} We assume that the following holds:
  \begin{itemize}
  \item {\bf (Independence):} All random variables in
    $\{\xi_{k,m}\}_{k=0,1,\cdots,K;m=1,\cdots, M}$ in Algorithm
    \ref{alg:conread} are independent to each other;
  \item {\bf (Bounded Age):} All delay variables $\tau_{k,m}$'s are
    bounded: $\max_{k,m} \tau_{k,m} \leq T$.
  \end{itemize}
\end{assumption}
The independence assumption strictly holds if all workers select
samples with \emph{replacement}. Although it might not be satisfied
strictly in practice, it is a common assumption made for the analysis
purpose. The bounded delay assumption is much more important. As
pointed out before, the asynchronous implementation may use some old
value of parameter $x$ to evaluate the stochastic
gradient. Intuitively, the age (or ``oldness'') should not be too
large to ensure the convergence. Therefore, it is a natural and
reasonable idea to assume an upper bound for ages. This assumption is
commonly used in the analysis for asynchronous algorithms, for
example, \citep{recht2011hogwild, avron2014revisiting,
  liu2014asynchronous, liu2013asynchronous, 2015arXiv150504824F,
  liu2014asynchronousRK}. It is worth noting that the upper bound $T$
is roughly proportional to the number of
workers.%

Under Assumptions~\ref{ass:0} and \ref{ass:1}, we have the following
convergence rate for nonconvex optimization.

\begin{theorem} \label{thm:conread} Assume that
  Assumptions~\ref{ass:0} and \ref{ass:1} hold and the steplength
  sequence $\{\gamma_k\}_{k=1,\cdots, K}$ in Algorithm
  \ref{alg:conread} satisfies
  \begin{align}
    LM\gamma_k+ 2L^2M^2T\gamma_k \sum_{\kappa=1}^{T}\gamma_{k+\kappa}  \leq 1 \quad \text{for all $k=1,2,...$.}
    \label{algo1ass}
  \end{align}
  We have the following ergodic convergence rate for the
  iteration of Algorithm~\ref{alg:conread}
  \begin{equation}
    \frac{1}{\sum_{k=1}^K\gamma_k}\sum_{k=1}^K\gamma_k \mathbb{E}(\left\| \nabla f(x_k) \right\|^2) \leq \frac{2(f(x_1) -f(x^*) ) + \sum_{k=1}^K\left(\gamma^2_k ML +  2L^2M^2\gamma_k \sum_{j=k-T}^{k-1}\gamma_j^2 \right)\sigma^2}{M\sum_{k=1}^K\gamma_k}.
    \label{eq:thm_1}
  \end{equation}
  where $\E(\cdot)$ denotes taking expectation in terms of all random
  variables in Algorithm \ref{alg:conread}.
\end{theorem}
To evaluate the convergence rate, the commonly used metrics in convex
optimization are not eligible, for example, $f(x_k)- f^*$
and $\|x_k - x^*\|^2$.
For nonsmooth optimization, we use the ergodic convergence as the
metric, that is, the weighted average of the $\ell_2$
norm of all gradients $\|\nabla f(x_k)\|^2$,
which is used in the analysis for nonconvex optimization
\citep{ghadimi2013stochastic}. Although the metric used in nonconvex
optimization is not exactly comparable to $f(x_k)- f^*$
or $\|x_k - x^*\|^2$
used in the analysis for convex optimization, it is not totally
unreasonable to think that they are roughly in the same order. The
ergodic convergence directly indicates the following convergence: If
randomly select an index $\tilde{K}$
from $\{1,2,\cdots, K\}$
with probability $\{\gamma_k/\sum_{k=1}^K \gamma_k\}$,
then $\E(\|\nabla f(x_{\tilde{K}})\|^2)$
is bounded by the right hand side of \eqref{eq:thm_1} and all bounds
we will show in the following.

Taking a close look at Theorem~\ref{thm:conread},
we can properly choose the steplength $\gamma_k$
as a constant value and obtain the following convergence rate:
\begin{corollary} \label{coro:conread}
  Assume that Assumptions~\ref{ass:0} and \ref{ass:1} hold. Set the steplength $\gamma_k$ to be a
  constant $\gamma$
  \begin{align}
    \gamma := \sqrt{f(x_1) - f(x^*) \over MLK\sigma^2}.
    \label{eq:coro_1}
  \end{align}
  If the delay parameter $T$ is bounded by
  \begin{align}
    K \geq \frac{4ML(f(x_1)- f(x^*))}{\sigma^2}(T+1)^2,
    \label{eq:coro_2}
  \end{align}
  then the output of Algorithm~\ref{alg:conread} satisfies the
  following ergodic convergence rate
  \begin{equation}
    \min_{k\in \{1, \cdots, K\}}\E(\|\nabla f(x_k)\|^2) \leq
    {1\over K}\sum_{k=1}^K \E(\|\nabla f(x_k)\|^2) \leq 4\sqrt{(f(x_1)- f(x^*))L \over MK}\sigma.
  \end{equation}
\end{corollary}
This corollary basically claims that when the total iteration number
$K$ is greater than $O(T^2)$, the convergence rate achieves
$O({1 / \sqrt{MK}})$. Since this rate does not depend on the delay
parameter $T$ after sufficient number of iterations, the negative
effect of using old values of $x$ for stochastic gradient evaluation
vanishes asymptoticly. In other words, if the total number of workers
is bounded by $O(\sqrt{K/M})$, the linear speedup is achieved.

Note that our convergence rate $O({1 / \sqrt{MK}})$
is consistent with the serial SG (with $M=1$)
for convex optimization \citep{nemirovski2009robust}, the synchronous
parallel (or mini-batch) SG for convex optimization
\citep{dekel2012optimal}, and nonconvex smooth optimization
\citep{ghadimi2013stochastic}. Therefore, an important observation is
that as long as the number of workers (which is
proportional to $T$)
is bounded by $O(\sqrt{K/M})$,
the iteration complexity to achieve the same accuracy level
will be roughly the same. In other words, the
average work load for each worker is reduced by the factor $T$
comparing to the serial SG. Therefore, the linear speedup is
achievable if $T\leq O(\sqrt{K/M})$.
Since our convergence rate meets several special cases, it is tight.

Next we compare with the analysis of $\AsySGC$ for \emph{convex} smooth
optimization in \citet[Corollary 2]{agarwal2011distributed}.
They proved an asymptotic convergence rate $O(1/\sqrt{MK})$,
which is consistent with ours. But their results require
$T\leq O(K^{1/4} M^{-3/4})$
to guarantee linear speedup. Our result improves it by a factor
$O(K^{1/4}M^{1/4})$.

\section{Asynchronous parallel stochastic gradient for shared memory
  architecture} \label{sec:asysgi}

This section considers a widely used lock-free asynchronous
implementation of SG on the shared memory system proposed in
\citet{recht2011hogwild}. Its advantages have been witnessed in solving
SVM, graph cuts \citep{recht2011hogwild}, linear equations
\citep{liu2014asynchronousRK}, and matrix completion
\citep{petroni2014gasgd}. While the computer network always involves
multiple machines, the shared memory platform usually only includes a
single machine with multiple cores / GPUs sharing the same memory.

\subsection{Algorithm Description: $\AsySGI$}

\begin{algorithm} [H]%
  \caption{$\AsySGI$} %
  \label{alg:inconread} %
  \begin{algorithmic} [1] %
    \REQUIRE $x_0$, $K$, $\gamma$ \ENSURE $x_K$
    \FOR{$k=0,\cdots, K-1$} \STATE Randomly select $M$ training
    samples indexed by $\xi_{k,1}, \xi_{k,2}, ... \xi_{k,M}$; \STATE
    Randomly select $i_k \in \{1,2,...,n\}$ with uniform distribution;
    \STATE
    {\small$ (x_{k+1})_{i_k} = (x_{k})_{i_k} -
      \gamma\sum\limits_{m=1}^{M}(G(\hat{x}_{k,m};
      \xi_{k,m}))_{i_k}$};
    \ENDFOR
  \end{algorithmic}
\end{algorithm}

For the shared memory platform, one can exactly
follow $\AsySGC$
on the computer network using software locks, which is
expensive\footnote{The time consumed by locks is roughly equal to
  the time of $10^4$
  floating-point computation. The additional cost for using locks is
  the waiting time during which multiple worker
  access the same memory address.}. Therefore, in practice the lock
free asynchronous parallel implementation of SG is preferred. This
section considers the same implementation as
\citet{recht2011hogwild}, but provides a more
precise algorithm description $\AsySGI$
than $\hogwild$ proposed in \citet{recht2011hogwild}.

In this lock free implementation, the shared memory stores the
parameter ``$x$'' and allows all workers reading and modifying
parameter $x$ simultaneously without using locks. All workers repeat
the following steps independently, concurrently, and simultaneously:
\begin{itemize}
\item {\bf (Read):} read the parameter from the
  shared memory to the local memory \emph{without software locks} (we
  use $\hat{x}$ to denote its value);
\item {\bf (Compute):} sample a training data $\xi$ and use $\hat{x}$
  to compute the stochastic gradient $G(\hat{x}; \xi)$ \emph{locally};
\item {\bf (Update):} update parameter $x$ in the shared memory
  \emph{without software locks}
  $x \leftarrow x - \gamma G(\hat{x}; \xi)$.
\end{itemize}
Since we do not use locks in both ``read'' and ``update'' steps, it
means that multiple workers may manipulate the shared memory
simultaneously. It causes the ``inconsistent read'' at the ``read''
step, that is, the value of $\hat{x}$
read from the shared memory might not be any state of $x$
in the shared memory at any time point. For example, at time $0$,
the original value of $x$
in the shared memory is a two dimensional vector $[a, b]$;
at time $1$,
worker $W$
is running the ``read'' step and first reads $a$
from the shared memory; at time $2$,
worker $W'$
updates the first component of $x$
in the shared memory from $a$
to $a'$;
at time $2$,
worker $W'$
updates the second component of $x$
in the shared memory from $b$
to $b'$;
at time $3$,
worker $W$
reads the value of the second component of $x$
in the shared memory as $b'$.
In this case, worker $W$
eventually obtains the value of $\hat{x}$
as $[a, b']$,
which is not a real state of $x$
in the shared memory at any time point. Recall that in $\AsySGC$
the parameter value obtained by any worker is guaranteed to be some
real value of parameter $x$ at some time
point. %

To precisely characterize this implementation and especially represent
$\hat{x}$, we monitor the value of parameter $x$ in the shared
memory. We define one \emph{iteration} as a modification on any
\emph{single} component of $x$ in the shared memory since the update
on a single component can be considered to be atomic on GPUs and DSPs
\citep{recht2011hogwild}. We use $x_{k}$ to denote the value of
parameter $x$ in the shared memory after $k$ iterations and
$\hat{x}_k$ to denote the value read from the shared memory and used for
computing stochastic gradient at the $k$th iteration. $\hat{x}_k$ can be
represented by $x_k$ with a few earlier updates missing
\begin{align}
  \hat{x}_k = x_k - \sum_{j\in J(k)} (x_{j+1}-x_j)
\end{align}
where $J(k)\subset \{k-1, k, \cdots, 0\}$ is a subset of index numbers
of previous iterations. This way is also used in analyzing
asynchronous parallel coordinate descent algorithms in
\citep{avron2014revisiting, liu2014asynchronous}. The $k$th update
happened in the shared memory can be described as
\[
  (x_{k+1})_{i_k} = (x_k)_{i_k} - \gamma (G(\hat{x}_{k}; \xi_k))_{i_k}
\]
where $\xi_k$ denotes the index of the selected data and $i_k$ denotes
the index of the component being updated at $k$th iteration. In the
original analysis for the $\hogwild$ implementation
\citep{recht2011hogwild}, $\hat{x}_{k}$ is assumed to be some earlier
state of $x$ in the shared memory (that is, the consistent read) for
simpler analysis, although it is not true in practice.

One more complication is to apply the mini-batch strategy like
before. Since the ``update'' step needs physical modification in the
shared memory, it is usually much more time consuming than both
``read'' and ``compute'' steps are. If many workers run the ``update''
step simultaneously, the memory contention will seriously harm the
performance. To reduce the risk of memory contention, a common trick
is to ask each worker to gather multiple (say $M$) stochastic
gradients and write the shared memory only once. That is, in each
cycle, run both ``update'' and ``compute'' steps for $M$ times before
you run the ``update'' step. Thus, the mini-batch updates happen in
the shared memory can be written as
\begin{align}
  (x_{k+1})_{i_k} = (x_k)_{i_k} - \gamma \sum_{m=1}^M (G(\hat{x}_{k, m}; \xi_{k,m}))_{i_k}
\end{align}
where $i_k$ denotes the coordinate index updated at the $k$th
iteration, and $G(\hat{x}_{k, m} ; \xi_{k,m})$ is the $m$th stochastic
gradient computed from the data sample indexed by $\xi_{k,m}$ and the
parameter value denoted by $\hat{x}_{k,m}$ at the $k$th
iteration. $\hat{x}_{k,m}$ can be expressed by:
\begin{align}
  \hat{x}_{k,m} = x_k - \sum_{j\in J(k,m)} (x_{j+1}-x_j)
\end{align}
where $J(k,m)\subset \{k-1, k, \cdots, 0\}$ is a subset of index
numbers of previous iterations. The algorithm is summarized in
Algorithm~\ref{alg:inconread} from the view of the shared memory.

\subsection{Analysis for $\AsySGI$}
To analyze the $\AsySGI$, we need to make a few assumptions similar to
\citet{recht2011hogwild, liu2014asynchronousRK, avron2014revisiting,
  liu2014asynchronous}.
\begin{assumption} \label{ass:3} We assume that the following holds
  for Algorithm \ref{alg:inconread}:
  \begin{itemize}
  \item {\bf (Independence):} All groups of variables $\{i_k, \{\xi_{k,m}\}_{m=1}^M\}$ at different iterations from $k=1$ to $K$ are independent to each other. %
  \item {\bf (Bounded Age):} Let $T$ be the global bound for delay:
    $J(k, m) \subset \{k-1, ... k-T\}, \quad \forall k, \forall m$,
    so $| J(k, m) | \leq T$.
  \end{itemize}
  \label{ass:2}
\end{assumption}
The independence assumption might not be true in practice, but it is
probably the best assumption one can make in order to analyze the
asynchronous parallel SG algorithm. This assumption was also used in
the analysis for $\hogwild$
\citep{recht2011hogwild} and asynchronous randomized Kaczmarz
algorithm \citep{liu2014asynchronousRK}. The bounded delay assumption
basically restricts the age of all missing components in
$\hat{x}_{k,m}$
($\forall m,~\forall k$).
The upper bound ``$T$''
here serves a similar purpose as in
Assumption~\ref{ass:1}. Thus we abuse this notation in this
section. The value of $T$
is proportional to the number of workers and does not depend on the
size of mini-batch $M$.
The bounded age assumption is used in the analysis for asynchronous
stochastic coordinate descent with ``inconsistent read''
\citep{avron2014revisiting, liu2014asynchronous}. Under
Assumptions~\ref{ass:0} and \ref{ass:2}, we have the following
results:
\begin{theorem}\label{thm:algo2}
  Assume that Assumptions~\ref{ass:0} and \ref{ass:2} hold
  and the constant steplength $\gamma$ satisfies
  \begin{equation}
    \frac{2M^2T L_T^2(\sqrt{n}+T-1)\gamma^2}{n^{3/2}} +2M\Lmax\gamma \le  1.
    \label{alg2ass}
  \end{equation}We have the following ergodic convergence rate for Algorithm~\ref{alg:inconread}
  \begin{equation}
    \frac{1}{K}\sum_{t=1}^K \mathbb{E}\left(  \|\nabla f (x_t)\|^2\right)\leq  \frac{2n}{KM\gamma} (f(x_1) - f(x^*))+ \frac{L_T^2TM\gamma^2}{2n} \sigma^2 + \Lmax\gamma \sigma^2.
    \label{eq:algo2}
  \end{equation}
\end{theorem} Taking a close look at
Theorem~\ref{thm:algo2}, we can choose the
steplength $\gamma$ properly and obtain the following error bound:
\begin{corollary} \label{coro:inconread}
  Assume that Assumptions~\ref{ass:0} and \ref{ass:2} hold. Set the steplength to be a
  constant $\gamma$
  \begin{align}
    \gamma := \frac{\sqrt{(f(x_1) - f(x^*))n}}{\sqrt{K L_T M} \sigma}.
    \label{eq:coro_11}
  \end{align}
  If the total iterations $K$ is greater than
  \begin{align}
    K \ge \frac{16 (f(x_1) - f(x^*)) L_T M \left(n^{3/2}+4 T^2\right)}{\sqrt{n} \sigma ^2},
    \label{eq:coro_22}
  \end{align}
  then the output of Algorithm~\ref{alg:inconread} satisfies the
  following ergodic convergence rate
  \begin{equation} {1\over K}\sum_{k=1}^K \E(\|\nabla f(x_k)\|^2) \leq
    \sqrt{\frac{72 \left(f\left(x_1\right) - f\left(x^*\right)\right) L_T n }{KM}}\sigma.
    \label{eq:coro_23}
  \end{equation}
\end{corollary}
This corollary indicates the asymptotic convergence rate achieves
$O(1/\sqrt{MK})$ when the total iteration number $K$ exceeds a
threshold in the order of $O(T^2)$ (if $n$ is considered as a
constant). We can see that this rate and the threshold are consistent
with the result in Corollary~\ref{coro:conread} for $\AsySGC$. One may
argue that why there is an additional factor $\sqrt{n}$ in the
numerator of \eqref{eq:coro_23}. That is due to the way we count
iterations -- one iteration is defined as updating a single component
of $x$. If we take into account this factor in the
comparison to $\AsySGC$, the convergence rates for $\AsySGC$ and
$\AsySGI$ are essentially consistent. This comparison implies that the
``inconsistent read'' would not make a big difference from the
``consistent read''.

Next we compare our result with the analysis of $\hogwild$ by
\citep{recht2011hogwild}. In principle, our analysis and their
analysis consider the same implementation of asynchronous parallel SG,
but differ in the following aspects: 1) our analysis considers the
smooth nonconvex optimization which includes the smooth strongly convex
optimization considered in their analysis;
2) our analysis considers the ``inconsistent read'' model which meets
the practice while their analysis assumes the impractical ``consistent
read'' model. Although the two results are not
absolutely comparable, it is still interesting to see the
difference. \citet{recht2011hogwild} proved that the linear speedup is
achievable if the maximal number of nonzeros in stochastic
gradients is bounded by $O(1)$ and the number of workers is bounded by
$O(n^{1/4})$. Our analysis does not need this prerequisite and
guarantees the linear speedup as long as the number of workers is
bounded by $O(\sqrt{K})$. Although it is hard to say that our result
strictly dominates $\hogwild$ in \citet{recht2011hogwild}, our
asymptotic result is eligible for more scenarios.

We consider an extension of $\AsySGI$
for sparse stochastic gradients. In this scenario,
we slightly change Steps 3 and 4 in Algorithm~\ref{alg:inconread}:
\\
\begin{center}
  \fbox{\parbox{0.8\textwidth}{3: Uniformly select $i_k$ from the support set of $g_k:=\sum_{m=1}^M G(\hat{x}_{k,m}; \xi_{k,m})$;\\
      ~~~~4:
      $(x_{k+1})_{i_k} = (x_{k})_{i_k} - \gamma\|g_k\|_0(g_k)_{i_k}$.
    }}
\end{center}
The convergence rate in \eqref{eq:coro_23} can be slightly improved by
taking the advantage of sparsity:
\begin{equation} {1\over K}\sum_{k=1}^K \E(\|\nabla f(x_k)\|^2) \leq
  \frac{6 \sqrt{2 \left(f\left(x_1\right) - f\left(x^*\right)\right) L_T \max_k\|g_k\|_0 } \sigma}{\sqrt{KM}}.
\end{equation}
The proof can be simply get by extending the proof
for \eqref{eq:coro_23}.

\section{Experiments}
\label{sec:exp}

The successes of $\AsySGC$ and $\AsySGI$ and their advantages over
synchronous parallel algorithms have been widely witnessed in many
applications such as deep neural network \citep{dean2012large,
  paine2013gpu, ZhangCL14a, li2014scaling}, matrix completion
\citep{recht2011hogwild, petroni2014gasgd, yun2013nomad}, SVM
\citep{recht2011hogwild}, and solving linear equations
\citep{liu2014asynchronousRK}. We refer readers to these literatures
for more comphrehensive comparison and empirical studies. This section
mainly provides the empirical study to validate the speedup properties
for \emph{completeness}.

We perform experiments for $\AsySGC$
and $\AsySGI$
on computer cluster and multicore machine respectively. The main
purpose of the following experiments is to validate the speedup
property. We are particularly interested in two types of speedup:
iteration speedup and running time speedup. The iteration speedup is
exactly the speedup we discussed in the whole paper. Given $T$
workers, it is computed from the ratio
\[
  \text{iteration speedup of $T$ workers} = \frac{\text{\# of total
      iterations of the serial SG (or using one worker)}}{\text{\# of
      total iterations using $T$ workers}} \times T
\]
where $\#$ is the iteration count when the same level of precision
achieved. This speedup is less affected by the hardware. The running
time speedup is the actual speedup. It is defined with respect to the
running time:
\[
  \text{running time speedup of $T$ workers} = \frac{\text{running
      time for the serial SG (or using one worker)}}{\text{running
      time of using $T$ workers}}.
\]
The running time speedup is seriously affected by the hardware. It is
generally worse than the iteration speedup.

\subsection{$\AsySGC$}

We implement $\AsySGC$
for deep neural network based on the Caffe \citep{jia2014caffe}
package. Caffe is an open source code base of deep learning
algorithms. We evaluate $\AsySGC$
on two standard datasets provided in the Caffe package, LENET and
CIFAR10-FULL.

The neural network consists of convolution layer, nonlinear layer, max
pool layer and fully connected layer \citep{NIPS2012_4824}, and the
detailed specification can be found on the Caffe
website\footnote{\url{https://github.com/BVLC/caffe}}. LENET is a
digit classifier network, training on the MNIST
dataset\footnote{\url{http://yann.lecun.com/exdb/mnist/}}.
CIFAR10-full has 10 classes of color images, training on the CIFAR10
dataset \citep{krizhevsky2009learning}. We first initialize a
parameter server hosting the parameters, and then spawn up to 8
stochastic gradient workers. The point to point communication between
the parameter server and gradient workers are handled by the MPICH
library\footnote{\url{https://www.mpich.org/}}. The parameter server
and stochastic gradient workers run on separate machines, and each
process uses a single core of a Xeon(R) E5-2430 CPU. The steplength
$\gamma$
for LENET is chosen as the default value. It means that this
steplength has been tuned to be the optimal for the serial SG
algorithm on LENET. The steplength we used for CIFAR10-FULL is chosen
as the default value as well.

\begin{figure} [htp!]
  \centering
  \includegraphics[page=1,width=0.48\textwidth]{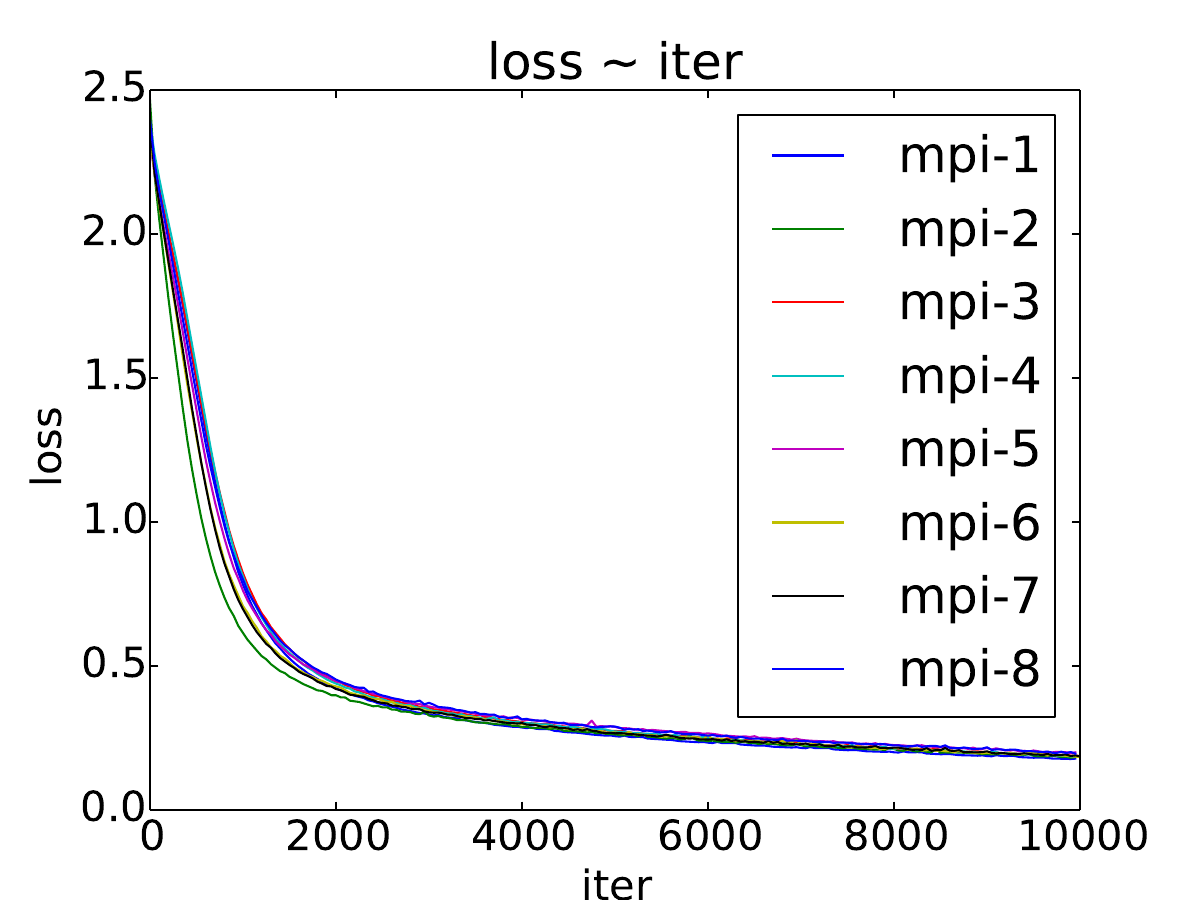}
  \includegraphics[page=2,width=0.48\textwidth]{figure/lenet-plots.pdf}
  \caption{LENET. The $\AsySGC$ algorithm is run on various numbers of machines from $1$ to $8$ to solve LENET. The curves of the objective loss against the number of iteration and the running time are drawn in the left and the right graphs respectively.} \label{fig:lenet}
  \includegraphics[page=1,width=0.48\textwidth]{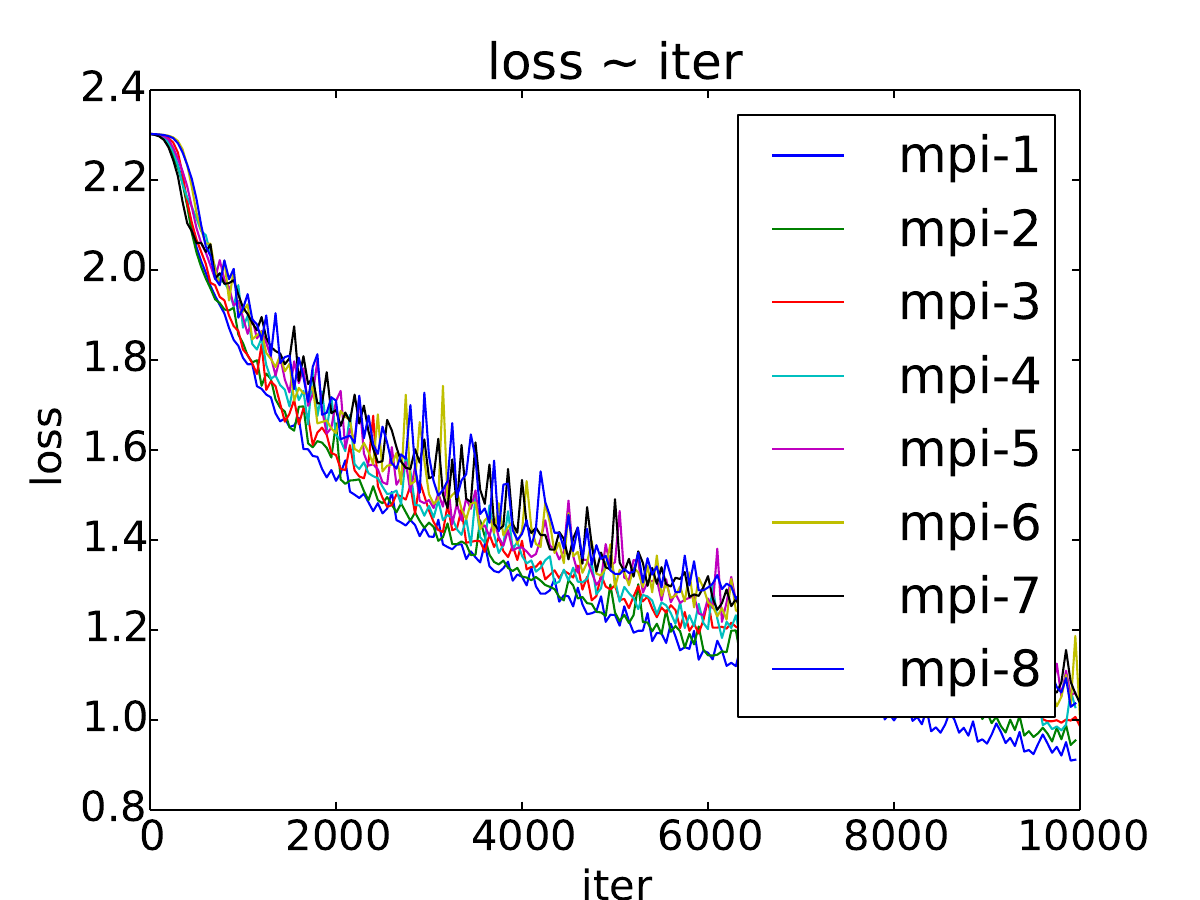}
  \includegraphics[page=2,width=0.48\textwidth]{figure/cifar10_full-plots.pdf}
  \caption{CIFAR10-FULL. The $\AsySGC$ algorithm is run on various numbers of machines from $1$ to $8$ to solve CIFAR10-FULL. The curves of the objective loss against the number of iteration and the running time are drawn in the left and the right graphs respectively.} \label{fig:cifar10-full}
\end{figure}

We draw the curves of objective loss against iterations and running
time in Figures \ref{fig:lenet} and \ref{fig:cifar10-full}
respectively, and report their speedups in Tables~\ref{tab:lenet} and
\ref{tab:cifar10-full}. (More details about datasets and the parameter
setting in experiments can be found in Table~\ref{tab:mpi-info}.) We
can observe that
\begin{itemize}
\item The iteration speedup is always better than the running time
  speedup, which meets our common sense;
\item The speedups for both problems are comparable overall as shown
  in Tables \ref{tab:lenet} and \ref{fig:cifar10-full}. The time
  speedup for CIFAR10-FULL is slightly more stable, while we notice
  that the time speedup for LENET in Table~\ref{tab:lenet} suddenly
  drops to ``$2.88$''
  with mpi-8 from ``5.29'' with mpi-7. That is because it hits the
  ceiling of communication bandwidth. The number of parameters in
  LENET is more than CIFAR10-FULL, thus requiring more communication
  cost. When the number of machines achieves a certain threshold (in
  this case LENET, it is $8$),
  the performance might become dramatically worse.
\end{itemize}

\begin{table} [htp!]
  \begin{center}
    \caption{Iteration speedup and running time speedup of
      $\AsySGC$. (LENET)} \label{tab:lenet}
    \begin{tabular}{ l c c c c c c c c }
      \toprule
      & mpi-1 & mpi-2 & mpi-3 & mpi-4 & mpi-5 & mpi-6 & mpi-7 & mpi-8  \\ \midrule
      iteration speedup & 1.07 & 2.02 & 2.77 & 3.73 & 4.20 & 5.82 & 6.86 & 6.97 \\ \hline
      time speedup & 0.98 & 1.69 & 2.24 & 2.96 & 3.30 & 4.51 & 5.29 & 2.88 \\
      \bottomrule
    \end{tabular}
  \end{center}
  \begin{center}
    \caption{Iteration speedup and running time speedup of
      $\AsySGC$. (CIFAR10-FULL)} \label{tab:cifar10-full}
    \begin{tabular}{ l c c c c c c c c }
      \toprule
      & mpi-1 & mpi-2 & mpi-3 & mpi-4 & mpi-5 & mpi-6 & mpi-7 & mpi-8  \\ \midrule
      iteration speedup & 1.01 & 1.93 & 2.65 & 3.42 & 4.27 & 4.92 & 5.36 & 5.96 \\ \hline
      time speedup & 1.00 & 1.73 & 2.28 & 2.88 & 3.56 & 4.07 & 4.41 & 5.00 \\
      \bottomrule
    \end{tabular}
  \end{center}
\end{table}

\begin{table} [htp!]  {%
    \centering
    \caption{More details about LENET and CIFAR10-FULL.} \label{tab:mpi-info}
    \begin{tabular}{c |cccccccc} \toprule
      Datasets & Type & \#Images & MiniBatch & \#CONV & \#FC & \#Params \\ \hline
      LENET & 28x28 grayscale & 60K & 60 & 2 & 2 & 431,080 \\
      CIFAR10-FULL & 32x32 RGB & 50K & 100 & 3 & 1 & 89,578 \\
      \bottomrule
    \end{tabular}
  }
\end{table}

\subsection{$\AsySGI$}
We conduct the empirical study for $\AsySGI$
on the machine (Intel Xeon architecture), which has 4 sockets and 10
cores for each socket. The synthetic data is generated from a full
connected neural network with 5 layers
($400\times 100 \times 50 \times 20 \times 10$)
and $46380$
parameters totally. The total number of samples is $463800$.
The data size is about $1.5$
GB. The input vector and all parameters are generated from
i.i.d. Gaussian distribution. The output vector is constructed by
applying the network parameter to the input vector plus some Gaussian
random noise.

\begin{figure}[htp!]
  \centering
  \includegraphics[width=0.48\textwidth]{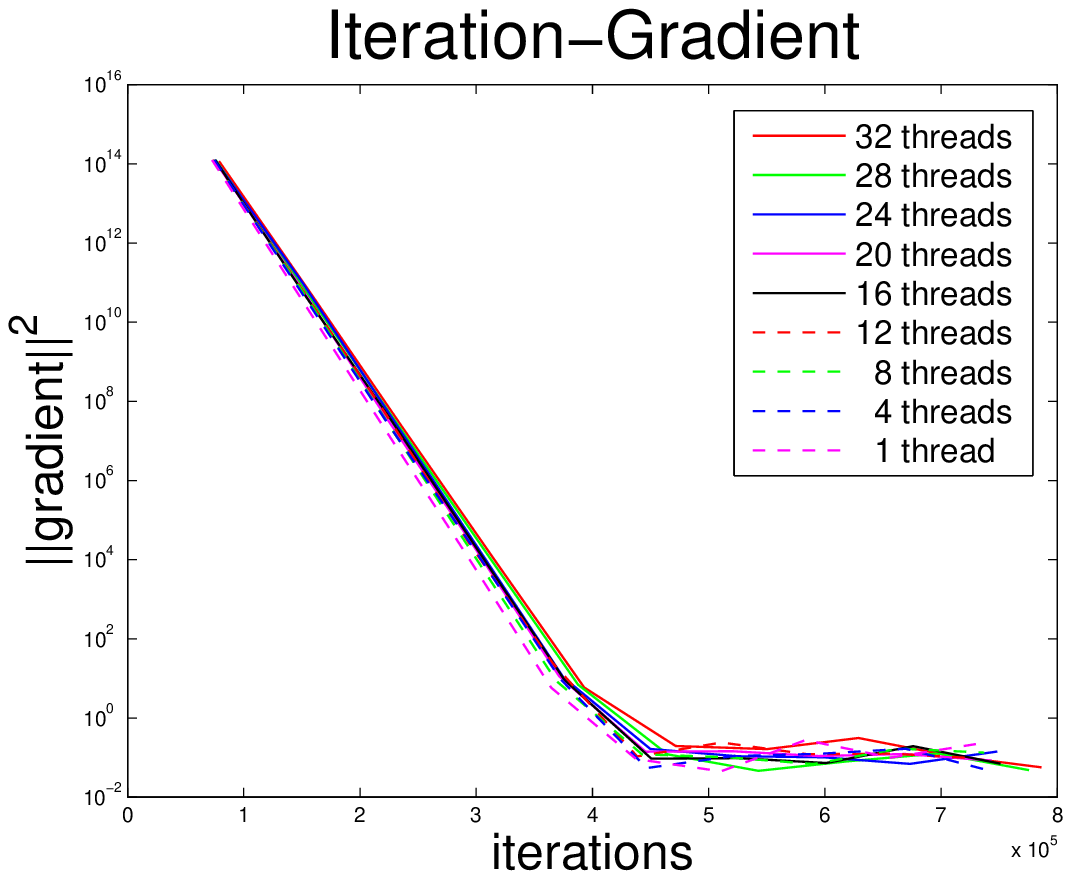}
  \includegraphics[width=0.48\textwidth]{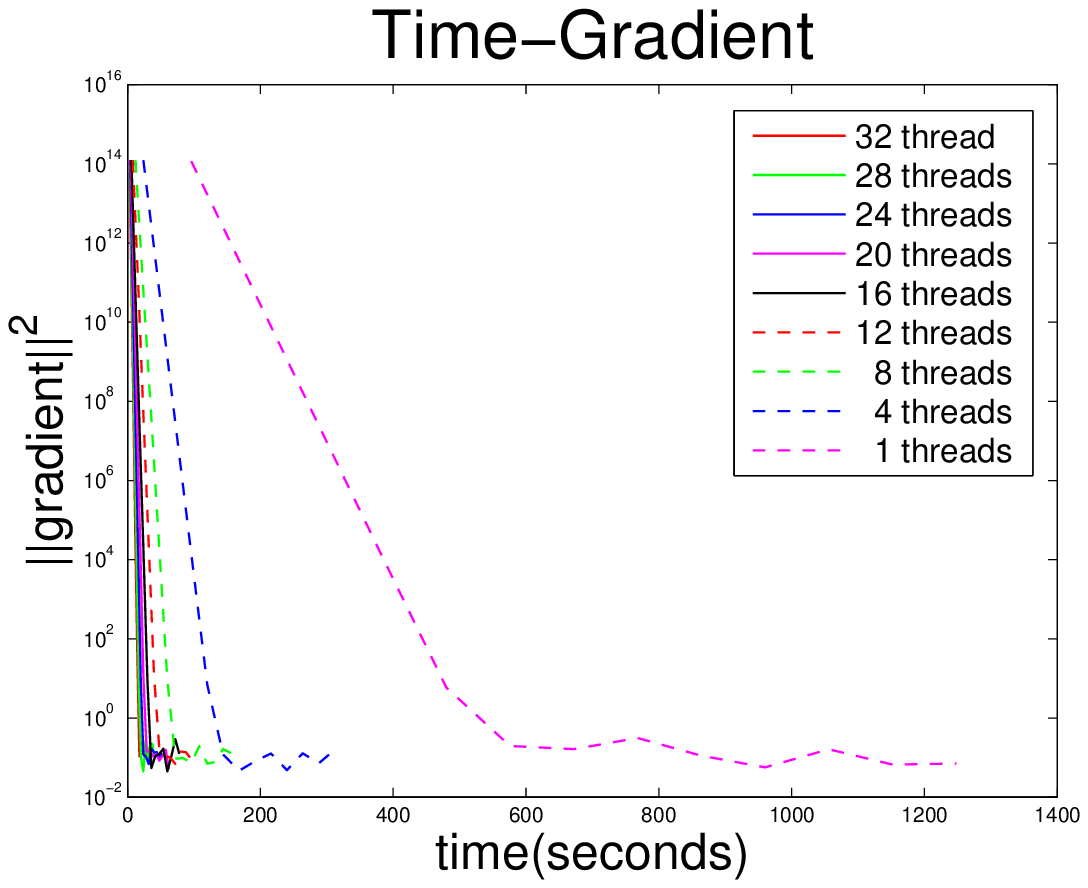}
  \caption{Deep neural network for synthetic data using $\AsySGI$.
    The $\AsySGC$
    algorithm is run on various numbers of machines from $1$
    to $32$.
    The curves of the objective loss against the number of iteration
    and the running time are drawn in the left and the right graphs
    respectively.}  \label{fig:asysgi:syn} \label{fig:asysgi:iter}
\end{figure}

We run $\AsySGI$
on various numbers of cores from $1$
to $32$.
The size of mini-batch is chosen as $M=32$
and the steplength is chosen as $\gamma=1.1\times 10^{-7}$.
Both parameters are chosen based on the best performance of the serial
SG to achieve the precision $10^{-1}$
for the $\ell_2$
norm of gradient. Figure~\ref{fig:asysgi:iter} draws the curve of the
$\ell_2$
norm of gradients against the number of iterations and running time
respectively. The speedup is reported in
Table~\ref{tab:asysgi:iter}. We observe that the iteration speedup is
almost linear while the running time speedup is slightly worse than
the iteration speedup. The overall performance of $\AsySGI$
on the shared memory system is better than $\AsySGC$
on the computer cluster. The reason is that the computer cluster
suffers from serious communication delay and the delay bound $T$
on the computer cluster is usually larger than on the shared memory
system.

\begin{table} [htp!]
  \begin{center}
    \caption{Iteration speedup and running time speedup of $\AsySGI$ (synthetic data).} \label{tab:asysgi:iter}
    \begin{tabular}{  l  c c  c  c  c  c  c  c  c }
      \toprule
      & thr-1 & thr-4 & thr-8 & thr-12 & thr-16 & thr-20 & thr-24 & thr-28 & thr-32  \\ \midrule
      iteration speedup &  1 & 3.9 & 7.8 & 11.6 & 15.4 & 19.9 & 24.1 & 28.7 & 31.6\\ \hline
      time speedup & 1 & 4.0 & 8.1 & 11.9 & 16.3 & 19.2 & 22.7 & 26.1 & 29.2\\
      \bottomrule
    \end{tabular}
  \end{center}
\end{table}

\section{Conclusion}
This paper studied two popular asynchronous parallel implementations
for $SG$
on computer cluster and shared memory system respectively. Two
algorithms ($\AsySGC$
and $\AsySGI$)
are used to describe two implementations. An asymptotic sublinear
convergence rate is proven for both algorithms on \emph{nonconvex}
smooth optimization. This rate is consistent with the result of $SG$
for convex optimization. The linear speedup is proven to achievable
when the number of workers is bounded by $\sqrt{K}$,
which improves the earlier analysis of $\AsySGC$
for convex optimization in \citep{agarwal2011distributed}. The
proposed $\AsySGI$
algorithm provides a more precise description for lock free
implementation on shared memory system than $\hogwild$
\citep{recht2011hogwild}. Our result for $\AsySGI$
can be applied to more scenarios.

\section*{Acknowledgements}
This project is supported by the NSF grant CNS-1548078, the NEC
fellowship, and the startup funding at University of Rochester. We
thank Professor Daniel Gildea and Professor Sandhya Dwarkadas at
University of Rochester, Professor Stephen J. Wright at University of
Wisconsin-Madison, and anonymous (meta-)reviewers for their constructive
comments and helpful advices.

\bibliographystyle{abbrvnat}

\newpage
\appendix
\section*{Appendix: Proofs}

{\noindent \bf Proofs to Theorem~\ref{thm:conread}}
\begin{proof}
  From the Lipschitzisan gradient assumption \eqref{ass:algo2-lip}, we have
  \begin{align}
    \nonumber
    f(x_{k+1}) -f(x_k)  \le &  \langle \nabla f(x_k), x_{k+1} - x_k \rangle + \frac{L}{2}\|x_{k+1} - x_k\|^2 \\
    = &
        - \left\langle \nabla f(x_k), \gamma_k \sum\limits_{m=1}^{M}G(x_{k-\tau_{k, m}}; \xi_{k, m})\right\rangle + \frac{\gamma_k^2 L }{2} \left\|\sum\limits_{m=1}^{M}G(x_{k-\tau_{k, m}}; \xi_{k, m})\right\|^2.
        \label{eq:proof:con_1}
  \end{align}
  Taking expectation respect to $\xi_{k, *}$ on both sides of
  \eqref{eq:proof:con_1}, we have
  \begin{align}
    \mathbb{E}_{\xi_{k, *}}(f(x_{k+1})) - f(x_k)
    \leq&  - M\gamma_k \left\langle \nabla f(x_k), \frac{1}{M}
          \sum_{m=1}^M\nabla f(x_{k-\tau_{k, m}}) \right\rangle \nonumber\\
        &+ \frac{\gamma^2_k L}{2}\mathbb{E}_{\xi_k, *}\left(\left\|\sum\limits_{m=1}^{M}G(x_{k-\tau_{k, m}}; \xi_{k, m})\right\|^2\right)
  \end{align}
  where we use the unbiased stochastic gradient assumption in \eqref{ass:1.2}. From the fact
  \[
\langle a ,b \rangle = \frac{1}{2}\left(\|a\|^2 + \|b\|^2 -
    \|a-b\|^2\right),\]
we have
  \begin{align}
    &\mathbb{E}_{\xi_{k, *}}(f(x_{k+1})) - f(x_k) \nonumber\\
    \leq& - \frac{M\gamma_k}{2} \left(\left\| \nabla f(x_k) \right\|^2 + \left\| \frac{1}{M} \sum_{m=1}^M\nabla f(x_{k-\tau_{k, m}}) \right\|^2 - \underbrace{\left\|\nabla f(x_k)-  \frac{1}{M} \sum_{m=1}^M\nabla f(x_{k-\tau_{k, m}}) \right\|^2}_{T_1}\right) \nonumber\\
    & + \frac{\gamma^2_k
      L}{2}\underbrace{\mathbb{E}_{\xi_k,*}\left(\left\|\sum\limits_{m=1}^{M}G(x_{k-\tau_{k,
      m}}; \xi_{k, m})\right\|^2\right)}_{T_2}.
      \label{eq:3}
  \end{align}

  Next we estimate the upper bound of $T_1$ and $T_2$. For $T_2$ we
  have
  \begin{align} T_2 =&\mathbb{E}_{\xi_{k, *}}\left[\left\|\sum\limits_{m=1}^{M}G\left(x_{{k}-\tau_{k, m}}; \xi_{k, m}\right)\right\|^2\right]\nonumber\\
    =& \mathbb{E}_{\xi_{k, *}} \left[ \left\| \sum_{m=1}^M \left(G(x_{{k}-\tau_{k, m}}; \xi_{k, m}) - \nabla f(x_{{k}-\tau_{k, m}})\right) + \sum_{m=1}^M \nabla f(x_{{k}-\tau_{k, m}}) \right\|^2\right]\nonumber\\
    =& \mathbb{E}_{\xi_{k, *}} \Bigg[ \left\| \sum_{m=1}^M \left(G(x_{{k}-\tau_{k, m}}; \xi_{k, m}) - \nabla f(x_{{k}-\tau_{k, m}})\right) \right\|^2 \nonumber\\
                     &+ \left\|\sum_{m=1}^M \nabla f(x_{{k}-\tau_{k,
                       m}}) \right\|^2 + 2 \left\langle \sum_{m=1}^M
                       \left(G(x_{{k}-\tau_{k, m}}; \xi_{k, m}) - \nabla f(x_{{k}-\tau_{k, m}})\right), \sum_{m=1}^M \nabla f(x_{{k}-\tau_{k, m}}) \right\rangle\Bigg]\nonumber\\
    =&  \mathbb{E}_{\xi_{k, *}} \left[ \left\| \sum_{m=1}^M \left(G(x_{{k}-\tau_{k, m}}; \xi_{k, m}) - \nabla f(x_{{k}-\tau_{k, m}})\right) \right\|^2 + \left\|\sum_{m=1}^M \nabla f(x_{{k}-\tau_{k, m}}) \right\|^2\right]\nonumber\\
    =& \mathbb{E}_{\xi_{k, *}} \Bigg[  \sum_{m=1}^M \left\|\left(G(x_{{k}-\tau_{k, m}}; \xi_{k, m})  - \nabla f(x_{{k}-\tau_{k, m}})\right) \right\|^2 \nonumber\\
                     &+ 2\sum_{1\le m < m^\prime \le M} \left\langle G(x_{{k}-\tau_{k, m}}; \xi_{k, m})  - \nabla f(x_{{k}-\tau_{k, m}}), G(x_{{k}-\tau_{k, m^{\prime}}}; \xi_{k, {m^\prime}})  - \nabla f(x_{{k}-\tau_{k, m^{\prime}}}) \right\rangle\nonumber\\
                     &+ \left\|\sum_{m=1}^M \nabla f(x_{{k}-\tau_{k, m}}) \right\|^2\Bigg]\nonumber\\
    \le &M \sigma^2 + \left\|\sum_{m=1}^M \nabla f(x_{{k}-\tau_{k,
          m}}) \right\|^2
          \label{eq:2}
  \end{align}
  where the forth equality is due to
  \begin{align*}
    &\mathbb{E}_{\xi_{k, *}} \left\langle \sum_{m=1}^M \left(G(x_{{k}-\tau_{k, m}}; \xi_{k, m}) - \nabla f(x_{{k}-\tau_{k, m}})\right), \sum_{m=1}^M \nabla f(x_{{k}-\tau_{k, m}}) \right \rangle \\
    =&  \left\langle \sum_{m=1}^M \mathbb{E}_{\xi_{k, *}}\left(G(x_{{k}-\tau_{k, m}}; \xi_{k, m}) - \nabla f(x_{{k}-\tau_{k, m}})\right), \sum_{m=1}^M \nabla f(x_{{k}-\tau_{k, m}}) \right \rangle \\
    =& 0
  \end{align*}
  and the last inequality is due to the assumption
  \eqref{eq:boundedvar} and
  \begin{align}
    &\mathbb{E}_{\xi_{k, *}}\sum_{1\le m < m^\prime \le M} \left\langle G(x_{{k}-\tau_{k, m}}; \xi_{k, m})  - \nabla f(x_{{k}-\tau_{k, m}}), G(x_{{k}-\tau_{k, m^{\prime}}}; \xi_{k, {m^\prime}})  - \nabla f(x_{{k}-\tau_{k, m^{\prime}}}) \right\rangle\nonumber\\
    =&  \mathbb{E}_{\xi_{k, *}}\sum_{1\le m < m^\prime \le M} \mathbb{E}_{k,m^\prime} \left\langle G(x_{{k}-\tau_{k, m}}; \xi_{k, m})  - \nabla f(x_{{k}-\tau_{k, m}}), G(x_{{k}-\tau_{k, m^{\prime}}}; \xi_{k, {m^\prime}})  - \nabla f(x_{{k}-\tau_{k, m^{\prime}}} \right\rangle\nonumber\\
    =& \mathbb{E}_{\xi_{k, *}}\sum_{1\le m < m^\prime \le M} \left\langle  G(x_{{k}-\tau_{k, m}}; \xi_{k, m})  - \nabla f(x_{{k}-\tau_{k, m}}), \mathbb{E}_{k,m^\prime} G(x_{{k}-\tau_{k, m^{\prime}}}; \xi_{k, {m^\prime}})  - \nabla f(x_{{k}-\tau_{k, m^{\prime}}}) \right\rangle\nonumber\\
    =& 0.
       \label{cross10}
  \end{align}

  We next turn to $T_1$:
  \begin{align*}
    T_1  =& \left\|\nabla f(x_k)- \frac{1}{M} \sum\limits_{m=1}^{M}\nabla f\left(x_{k-\tau_{k, m}}\right)\right\|^2 \\
    =& \frac{1}{M^2} \left\|\sum\limits_{m=1}^{M}\left(\nabla f(x_k)- \nabla f\left(x_{k-\tau_{k, m}}\right)\right)\right\|^2\\
    \leq& \frac{1}{M} \sum\limits_{m=1}^{M}\left\|\left(\nabla f(x_k)- \nabla f\left(x_{k-\tau_{k, m}}\right)\right)\right\|^2\\
    \leq & \frac{L^2}{M}\sum_{m=1}^M \|x_k - x_{{k}-\tau_{k, m}}\|^2 \\
    \leq&L^2\max_{k\in \{1,\cdots, M\}}\|x_k - x_{k-\tau_{k, m}}\|^2 \\
    =&L^2\|x_k - x_{k-\tau_{k, \mu}}\|^2. \quad \text{(let $\mu:={\arg\max}_{m\in \{1,\cdots, M\}}\|x_k - x_{k-\tau_{k, m}}\|^2$)}
  \end{align*}
  where the second inequality is from the Lipschitzian gradient
  assumption \eqref{ass:algo2-lip}. It follows that
  \begin{align}
    T_1 \leq& L^2\left\|x_k - x_{k-\tau_{k, \mu}}\right\|^2\nonumber\\
    =& L^2 \left\| \sum_{j=k-\tau_{k, \mu}}^{k-1} \left(x_{j+1} - x_j\right) \right\|^2\nonumber\\
    =& L^2 \left\|\sum_{j=k-\tau_{k, \mu}}^{k-1} \gamma_j\sum\limits_{m=1}^{M}G \left(x_{j- \tau_{j,m}}; \xi_{j,m}\right)\right\|^2 \nonumber\\
    =& L^2  \left\|\sum_{j=k-\tau_{k, \mu}}^{k-1}\gamma_j \sum\limits_{m=1}^{M} \left[G \left(x_{j- \tau_{j,m}}; \xi_{j,m}\right)- \nabla f\left(x_{j- \tau_{j,m}}\right) \right] + \sum_{j=k-\tau_{k, \mu}}^{k-1}  \gamma_j\sum\limits_{m=1}^{M}\nabla f\left(x_{j- \tau_{j,m}}\right) \right\|^2\nonumber \\
    \leq &2L^2  \left( \underbrace{\left\|\sum_{j=k-\tau_{k, \mu}}^{k-1}\gamma_j \sum\limits_{m=1}^{M} \left[G \left(x_{j- \tau_{j,m}}; \xi_{j,m}\right)- \nabla f\left(x_{j- \tau_{j,m}}\right) \right]\right\|^2}_{T_3} + \underbrace{\left\|\sum_{j=k-\tau_{k, \mu}}^{k-1} \gamma_j \sum\limits_{m=1}^{M}\nabla f\left(x_{j- \tau_{j,m}}\right) \right\|^2}_{T_4} \right)
           \label{eq:t1}
  \end{align}
  where the last inequality uses the fact that $\|a + b\|^2 \le 2\|a\|^2 + 2\|b\|^2$ for any real vectors $a$ and $b$. Taking the expectation in terms of $\{\xi_{j,*}|j\in \{k-\tau_{k, \mu}, ..., k-1\}\}$ for $T_3$, we have
  \begin{align}
    \nonumber    &\mathbb{E}_{\xi_{j,*},j\in \{k-\tau_{k, \mu},..., k-1\}}(T_3)
    \\      =&  \mathbb{E}_{\xi_{j,*},j\in \{k-\tau_{k, \mu},..., k-1\}}\left( \left\|\sum_{j=k-\tau_{k, \mu}}^{k-1}\gamma_j\sum\limits_{m=1}^{M} \left[G (x_{j - \tau_{j,m}}; \xi_{j,m}) - \nabla f(x_{j - \tau_{j,m}})\right ]\right\|^2\right)
               \nonumber \\
    =& \mathbb{E}_{\xi_{j,*},j\in \{k-\tau_{k, \mu},..., k-1\}}\left(  \sum_{j=k-\tau_{k, \mu}}^{k-1}\gamma_j^2 \left\|\sum\limits_{m=1}^{M}\left[G (x_{j - \tau_{j,m}}; \xi_{j,m}) - \nabla f(x_{j - \tau_{j,m}})\right ]\right\|^2 \right)
       \nonumber \\
                 &+2\mathbb{E}_{\xi_{j,*},j\in \{k-\tau_{k, \mu},..., k-1\}}\Bigg(\nonumber
                   \sum_{k-1\ge
                   j^{\prime\prime}
                   >
                   j^\prime\ge
                   k-\tau_{k,\mu}}\gamma_{j^{\prime}}\gamma_{j^{\prime\prime}} \Bigg\langle \sum\limits_{m=1}^{M}\left[G (x_{j^{\prime\prime} - \tau_{j^{\prime\prime},m}}; \xi_{j^{\prime\prime},m}) - \nabla f(x_{j^{\prime\prime} - \tau_{j^{\prime\prime},m}})\right ],\nonumber\\
                 &\sum\limits_{m=1}^{M}\left[G (x_{j^\prime -
                   \tau_{j^\prime,m}}; \xi_{j^\prime,m}) -
                   \nabla_{i_{j^\prime}}f(x_{j^{\prime}- \tau_{j^{\prime},m}})\right ]\Bigg\rangle \Bigg)
                   \nonumber\\
    =& \mathbb{E}_{\xi_{j,*},j\in \{k-\tau_{k, \mu},..., k-1\}}\left(  \sum_{j=k-\tau_{k, \mu}}^{k-1}\gamma_j^2 \left\|\sum\limits_{m=1}^{M}\left[G (x_{j - \tau_{j,m}}; \xi_{j,m}) - \nabla f(x_{j - \tau_{j,m}})\right ]\right\|^2 \right)  \nonumber\\
    =& \mathbb{E}_{\xi_{j,*},j\in \{k-\tau_{k, \mu},..., k-1\}}\left(
       \sum_{j=k-\tau_{k, \mu}}^{k-1}\gamma_j^2
       \sum\limits_{m=1}^{M}\left\|\left[G (x_{j - \tau_{j,m}};
       \xi_{j,m}) - \nabla f(x_{j - \tau_{j,m}})\right ]\right\|^2
       \right)
       \nonumber\\
    \leq & M\sum_{j=k-T}^{k-1}\gamma_j^2\sigma^2
           \label{eq:t3}
  \end{align}
  where the second last equality is due to \eqref{cross10} and the
  third equality is due to
  \begin{align*}
    &\mathbb{E}_{\xi_j,k-1\ge
      j
      \ge
      k-\tau_{k,\mu}}\Bigg(\sum_{k-1\ge
      j^{\prime\prime}
      >
      j^\prime\ge
      k-\tau_{k,\mu}}\gamma_{j^{\prime}}\gamma_{j^{\prime\prime}} \Bigg\langle \sum\limits_{m=1}^{M}\left[G (x_{j^{\prime\prime} - \tau_{j^{\prime\prime},m}}; \xi_{j^{\prime\prime},m}) - \nabla f(x_{j^{\prime\prime} - \tau_{j^{\prime\prime},m}})\right ],\nonumber\\
    &\sum\limits_{m=1}^{M}\left[G (x_{j^\prime - \tau_{j^\prime,m}};
      \xi_{j^\prime,m}) - \nabla f(x_{j^{\prime}- \tau_{j^{\prime},m}})\right ]\Bigg\rangle \Bigg)\\
    =&\mathbb{E}_{\xi_j,k-1\ge
       j
       \ge
       k-\tau_{k,\mu}}\Bigg(\sum_{k-1\ge
       j^{\prime\prime}
       >
       j^\prime\ge
       k-\tau_{k,\mu}}\gamma_{j^{\prime}}\gamma_{j^{\prime\prime}} \Bigg\langle \sum\limits_{m=1}^{M}\left[\mathbb{E}_{j^{\prime\prime}*}G (x_{j^{\prime\prime} - \tau_{j^{\prime\prime},m}}; \xi_{j^{\prime\prime},m}) - \nabla f(x_{j^{\prime\prime} - \tau_{j^{\prime\prime},m}})\right ],\nonumber\\
    &\sum\limits_{m=1}^{M}\left[G (x_{j^\prime - \tau_{j^\prime,m}};
      \xi_{j^\prime,m}) - \nabla f(x_{j^{\prime}- \tau_{j^{\prime},m}})\right ]\Bigg\rangle \Bigg)\\
    =&0.
  \end{align*}
  Taking the expectation in terms of $\xi_{j,*}$ for $T_4$, we have
  \begin{align}
    &\mathbb{E}_{\xi_{j,*},j\in \{k-\tau_{k, \mu},..., k-1\}}(T_4) \nonumber\\
    =&\mathbb{E}_{\xi_{j,*},j\in \{k-\tau_{k, \mu},..., k-1\}}\left( \left\|\sum_{j=k-\tau_{k, \mu}}^{k-1} \gamma_j \sum\limits_{m=1}^{M}\nabla f\left(x_{j- \tau_{j,m}}\right) \right\|^2 \right)
       \nonumber\\
    \leq &
           T\sum_{j=k-\tau_{k, \mu}}^{k-1} \gamma_j^2 \mathbb{E}_{\xi_{j,*},j\in \{k-\tau_{k, \mu},..., k-1\}}\left(\left\|\sum\limits_{m=1}^{M}\nabla f(x_{j - \tau_{j,m}})\right\|^2\right)
           \label{eq:t4}
  \end{align}
  where the last inequality uses the upper bound of the delay age:
  $\tau_{k,\mu} \le T$.

  We take full expectation on both sides of \eqref{eq:t1} and
  substitute $\mathbb{E}(T_{3})$ and $\mathbb{E}(T_{4})$ by their
  upper bounds in \eqref{eq:t3} and \eqref{eq:t4} respectively:

  \begin{align}
    \mathbb{E}(T_1) \leq 2L^2  \left( M\sum_{j=k-T}^{k-1}\gamma_j^2\sigma^2 + T\sum_{j=k-\tau_{k, \mu}}^{k-1} \gamma_j^2 \mathbb{E}\left(\left\|\sum\limits_{m=1}^{M}\nabla f(x_{j - \tau_{j,m}})\right\|^2\right)
    \right).
    \label{eq:1}
  \end{align}
  Applying the upper bounds for $\mathbb{E}(T_1)$ in \eqref{eq:1} and
  $\mathbb{E}(T_2)$ in \eqref{eq:2} to \eqref{eq:3}, and take full
  expectation on both sides, we obtain
  \begin{align}
    \nonumber&\mathbb{E}(f(x_{k+1})) - f(x_k)\\
    \nonumber\leq& - \frac{M\gamma_k}{2} \Bigg[\mathbb{E}\left(  \left\| \nabla f (x_k) \right\|^2\right) + \mathbb{E}\left(  \left\| \frac{1} {M} \sum_ {m=1}^ M\nabla f (x_{k-\tau_{k, m}}) \right\|^2\right)\\
    \nonumber& -2L^2  \left( M\sum_{j=k-T}^{k-1}\gamma_j^2\sigma^2 + T\sum_{j=k-\tau_{k, \mu}}^{k-1} \gamma_j^2 \mathbb{E}\left(\left\|\sum\limits_{m=1}^{M}\nabla f(x_{j - \tau_{j,m}})\right\|^2\right)
               \right)\Bigg] \\
    \nonumber& + \frac{\gamma^2_k L}{2}\left(M \sigma^2 +  \mathbb{E}\left(  \left\|\sum_ {m=1}^ M \nabla f (x_{{k}-\tau_{k, m}}) \right\|^2\right)\right)\\
    \nonumber\le & - \frac{M\gamma_k}{2}\mathbb{E}\left\| \nabla f(x_k) \right\|^2 + \left(\frac{\gamma^2_k L}{2}- \frac{\gamma_k}{2M}\right) \mathbb{E}\left(  \left\| \sum_{m=1}^ M\nabla f (x_{k-\tau_{k, m}}) \right\|^2\right)\\
    \nonumber&+ \left(\frac{\gamma^2_k ML}{2} + L^2M^2\gamma_k \sum_{j=k-T}^{k-1}\gamma_j^2 \right)\sigma^2\\
             &+L^2MT\gamma_k \sum_{j=k-T}^{k-1} \gamma_j^2 \mathbb{E}\left(\left\|\sum\limits_{m=1}^{M}\nabla f(x_{j - \tau_{j,m}})\right\|^2\right).
               \label{eq:proof:con:4}
  \end{align}
  Summarizing the inequality \eqref{eq:proof:con:4} from $k=1$ to $k=K$, we have
  \begin{align*}
    \nonumber&\mathbb{E}(f(x_{K+1})) - f(x_1)\\
    \nonumber\le & - \frac{M}{2}\sum_{k=1}^K\gamma_k\mathbb{E}\left(  \left\| \nabla f (x_k) \right\|^2\right) + \sum_{k=1}^K \left(\frac{\gamma^2_k L}{2}- \frac{\gamma_k}{2M}\right) \mathbb{E}\left(  \left\| \sum_ {m=1}^ M\nabla f (x_{k-\tau_{k, m}}) \right\|^2\right)\\
    \nonumber&+ \sum_{k=1}^K\left(\frac{\gamma^2_k ML}{2} + L^2M^2\gamma_k \sum_{j=k-T}^{k-1}\gamma_j^2 \right)\sigma^2\\
             &+L^2MT\sum_{k=1}^K\left(\gamma_k \sum_{j=k-T}^{k-1} \gamma_j^2 \mathbb{E}\left(\left\|\sum\limits_{m=1}^{M}\nabla f(x_{j - \tau_{j,m}})\right\|^2\right)\right)\\
    = & - \frac{M}{2}\sum_{k=1}^K\gamma_k\mathbb{E}\left(  \left\| \nabla f (x_k) \right\|^2\right)\\
             & + \sum_{k=1}^K \left(\gamma^2_k \left(\frac{L}{2} + L^2MT \sum_{\kappa=1}^{T}\gamma_{k+\kappa} \right)- \frac{\gamma_k}{2M}\right) \mathbb{E}\left(  \left\| \sum_{m=1}^M\nabla f (x_{k-\tau_{k, m}}) \right\|^2\right)\\
             & +\sum_{k=1}^K\left(\frac{\gamma^2_k ML}{2} + L^2M^2\gamma_k \sum_{j=k-T}^{k-1}\gamma_j^2 \right)\sigma^2\\
    \le& - \frac{M}{2}\sum_{k=1}^K\gamma_k\mathbb{E}\left(  \left\| \nabla f (x_k) \right\|^2\right) +\sum_{k=1}^K\left(\frac{\gamma^2_k ML}{2} + L^2M^2\gamma_k \sum_{j=k-T}^{k-1}\gamma_j^2 \right)\sigma^2
  \end{align*}
  where the last inequality is due to \eqref{algo1ass}. Note that
  $x^*$ is the global optimization point. Thus we have
  \begin{align*}
    \frac{1}{\sum_{k=1}^K\gamma_k}\sum_{k=1}^K\gamma_k \mathbb{E}(\left\| \nabla f(x_k) \right\|^2) \leq \frac{2(f(x_1) -f(x^*) ) + \sum_{k=1}^K\left(\gamma^2_k ML +  2L^2M^2\gamma_k \sum_{j=k-T}^{k-1}\gamma_j^2 \right)\sigma^2}{M\sum_{k=1}^K\gamma_k}.
  \end{align*}
  It completes the proof.
\end{proof}

{\noindent \bf Proofs to Corollary~\ref{coro:conread}}
\begin{proof}
  From \eqref{eq:coro_1} and \eqref{eq:coro_2}, we have
  \begin{equation}
    \gamma \leq {1\over 2ML(T+2)}.
    \label{eq:proof_con_1}
  \end{equation}
  If follows that
  \[ {1\over 2}\gamma L + L^2M T^2 \gamma^2 \le {1\over 4M (T+2)} +
    {T^2 \over 4M(T+2)^2}\leq {1\over 2M},
  \]
  which implies that the condition \eqref{algo1ass} in
  Theorem~\ref{thm:conread} is satisfied globally. Then we can safely
  apply \eqref{eq:thm_1} in Theorem~\ref{thm:conread}:
  \begin{align}
    \nonumber
    {1\over K}\sum_{i=1}^K \E(\|\nabla f(x_i)\|^2) \leq &\frac{2(f(x_1) -f(x^*) ) + K\left(\gamma^2 ML +  2L^2M^2T\gamma^3 \right)\sigma^2}{MK\gamma}
    \\ \nonumber = &
                     \frac{2(f(x_1)-f(x^*))}{MK\gamma} + L\sigma^2\gamma + 2L^2MT\sigma^2\gamma^2
    \\ \nonumber \leq &
                        \frac{2(f(x_1)- f(x^*))}{MK\gamma} + 2L\sigma^2\gamma
    \\ \nonumber = &
                     4\sqrt{(f(x_1)- f(x^*))L \over MK}\sigma,
  \end{align}
  where the second last inequality is due to \eqref{eq:proof_con_1}
  and the last equality uses \eqref{eq:coro_1}. It completes the
  proof.
\end{proof}

{\noindent \bf Proofs to Theorem~\ref{thm:algo2}}
\begin{proof}
  From the Lipschitzisan gradient assumption \eqref{ass:algo2-lip}, we
  have
  \begin{align}
    f(x_{k+1}) \leq f(x_k) - \gamma\left\langle\nabla_{i_k}f(x_k), \sum\limits_{m=1}^{M}(G(\hat{x}_{k, m} ; \xi_{k, m}))_{i_k}\right\rangle + \frac{\Lmax\gamma^2}{2}\left(\sum\limits_{m=1}^{M}(G(\hat{x}_{k, m} ; \xi_{k, m}))_{i_k}\right)^2.
    \label{eq:proof:incon:1}
  \end{align}
  Taking the expectation of $i_k$ on both sides of
  \eqref{eq:proof:incon:1}, we have
  \begin{align}
    \mathbb{E}_{i_k}f(x_{k+1}) \leq f(x_k) - {\gamma\over n}\left\langle\nabla f(x_k), \sum\limits_{m=1}^{M}G(\hat{x}_{k, m} ; \xi_{k, m})\right\rangle + \frac{\Lmax\gamma^2}{2n}\left\|\sum\limits_{m=1}^{M}G(\hat{x}_{k, m} ; \xi_{k, m})\right\|^2.
    \label{first_algo2}
  \end{align}
  Taking the expectation of $\xi_{k, *}$ on both sides of
  \eqref{first_algo2}, we obtain
  \begin{align}
    &\mathbb{E}_{\xi_{k, *}, i_k} \left(f\left(x_{k+1}\right)\right)
      \nonumber \\
    \leq& f\left(x_k\right) - \frac{\gamma}{n}\mathbb{E}_{\xi_{k, *}}\left[\left\langle\nabla f\left(x_k\right), \sum\limits_{m=1}^{M}G(\hat{x}_{k, m} ; \xi_{k, m})\right\rangle \right] +\frac{\Lmax\gamma^2}{2n}\mathbb{E}_{\xi_{k, *}} \left[\left\|\sum\limits_{m=1}^{M}G(\hat{x}_{k, m} ; \xi_{k, m})\right\|^2\right]\nonumber\\
    =&  f\left(x_k\right) - \frac{\gamma}{n} \left\langle \nabla f\left(x_k\right),\sum\limits_{m=1}^{M} \nabla f\left(\hat{x}_{k, m}\right) \right\rangle +\frac{\Lmax\gamma^2}{2n}\mathbb{E}_{\xi_{k, *}} \left[\left\|\sum\limits_{m=1}^{M}G(\hat{x}_{k, m} ; \xi_{k, m})\right\|^2\right]\nonumber\\
    =&f\left(x_k\right) - \frac{M\gamma}{2n}  \left(\left\|\nabla f\left(x_k\right)\right\|^2 + \left\|\frac{1}{M}\sum\limits_{m=1}^{M}\nabla f\left(\hat{x}_{k, m}\right)\right\|^2 - \underbrace{\left\|\nabla f\left(x_k\right) - \frac{1}{M} \sum\limits_{m=1}^{M}\nabla f\left(\hat{x}_{k, m}\right)\right\|^2}_{T_1}  \right)\nonumber\\
    &+\frac{\Lmax\gamma^2}{2n}\mathbb{E}_{\xi_{k, *}}\underbrace{ \left[\left\|\sum\limits_{m=1}^{M}G(\hat{x}_{k, m} ; \xi_{k, m})\right\|^2\right]}_{T_2}
      \label{eq:lip}
  \end{align}
  where the last equation is deduced by the fact
  $2\langle a, b\rangle = \|a\|^2+\|b\|^2-\|a-b\|^2$. We next consider
  $T_1$ and $T_2$ respectively.

  For $T_2$, we have
  \begin{align} \mathbb{E}_{\xi_{k, *}}\left(T_2\right)  =&\mathbb{E}_{\xi_{k, *}}\left[\left\|\sum\limits_{m=1}^{M}G(\hat{x}_{k, m} ; \xi_{k, m})\right\|^2\right]\nonumber\\
    =& \mathbb{E}_{\xi_{k, *}} \left[ \left\| \sum_{m=1}^M \left(G(\hat{x}_{k, m} ; \xi_{k, m}) - \nabla f(\hat{x}_{k, m})\right) + \sum_{m=1}^M \nabla f(\hat{x}_{k, m}) \right\|^2\right]\nonumber\\
    =& \mathbb{E}_{\xi_{k, *}} \Bigg[ \left\| \sum_{m=1}^M \left(G(\hat{x}_{k, m} ; \xi_{k, m}) - \nabla f(\hat{x}_{k, m})\right) \right\|^2 \nonumber\\
                                                          &+ \left\|\sum_{m=1}^M \nabla f(\hat{x}_{k, m}) \right\|^2 + 2 \left\langle \sum_{m=1}^M \left(G(\hat{x}_{k, m} ; \xi_{k, m}) - \nabla f(\hat{x}_{k, m})\right), \sum_{m=1}^M \nabla f(\hat{x}_{k, m}) \right\rangle\Bigg]\nonumber\\
    =&  \mathbb{E}_{\xi_{k, *}} \left[ \left\| \sum_{m=1}^M \left(G(\hat{x}_{k, m} ; \xi_{k, m}) - \nabla f(\hat{x}_{k, m})\right) \right\|^2 + \left\|\sum_{m=1}^M \nabla f(\hat{x}_{k, m}) \right\|^2\right]\nonumber\\
    =& \mathbb{E}_{\xi_{k, *}} \Bigg[  \sum_{m=1}^M \left\|G(\hat{x}_{k, m} ; \xi_{k, m})  - \nabla f(\hat{x}_{k, m}) \right\|^2 \nonumber\\
                                                          &+ 2\sum_{1\le m < m^\prime \le M} \left\langle G(\hat{x}_{k, m} ; \xi_{k, m})  - \nabla f(\hat{x}_{k, m}), G(\hat{x}_{k, {m^\prime}} ; \xi_{k, {m^\prime}})  - \nabla f(\hat{x}_{k, {m^\prime}}) \right\rangle
                                                            \nonumber   \\  & + \left\|\sum_{m=1}^M \nabla f(\hat{x}_{k, m}) \right\|^2\Bigg]\nonumber\\
    \le &M \sigma^2 + \left\|\sum_{m=1}^M \nabla f(\hat{x}_{k, m})
          \right\|^2
          \label{eq:10}
  \end{align}
  where the forth equality is due to
  \begin{align*}
    &\mathbb{E}_{\xi_{k, *}} \left\langle \sum_{m=1}^M \left(G(\hat{x}_{k, m} ; \xi_{k, m}) - \nabla f(\hat{x}_{k, m})\right), \sum_{m=1}^M \nabla f(\hat{x}_{k, m}) \right \rangle \\
    =&  \left\langle \sum_{m=1}^M \mathbb{E}_{\xi_{k, *}}\left(G(\hat{x}_{k, m} ; \xi_{k, m}) - \nabla f(\hat{x}_{k, m})\right), \sum_{m=1}^M \nabla f(\hat{x}_{k, m}) \right \rangle \\
    =& 0
  \end{align*}
  and the last inequality is due to the assumption
  \eqref{eq:boundedvar} and
  \begin{align}
    &\mathbb{E}_{\xi_{k, *}}\sum_{1\le m < m^\prime \le M} \left\langle G(\hat{x}_{k, m} ; \xi_{k, m})  - \nabla f(\hat{x}_{k, m}), G(\hat{x}_{k, {m^\prime}} ; \xi_{k, {m^\prime}})  - \nabla f(\hat{x}_{k, {m^\prime}}) \right\rangle\nonumber\\
    =&  \mathbb{E}_{\xi_{k, *}}\sum_{1\le m < m^\prime \le M} \mathbb{E}_{k, {m^\prime}} \left\langle G(\hat{x}_{k, m} ; \xi_{k, m})  - \nabla f(\hat{x}_{k, m}), G(\hat{x}_{k, {m^\prime}} ; \xi_{k, {m^\prime}})  - \nabla f(\hat{x}_{k, {m^\prime}}) \right\rangle\nonumber\\
    =& \mathbb{E}_{\xi_{k, *}}\sum_{1\le m < m^\prime \le M} \left\langle  G(\hat{x}_{k, m} ; \xi_{k, m})  - \nabla f(\hat{x}_{k, m}),  \mathbb{E}_{k, {m^\prime}}G(\hat{x}_{k, {m^\prime}} ; \xi_{k, {m^\prime}})  - \nabla f(\hat{x}_{k, {m^\prime}}) \right\rangle\nonumber\\
    =& 0.
       \label{cross0}
  \end{align}

  As for $T_1$, we have:
  \begin{align*}
    T_1  =& \left\|\nabla f(x_k)- \frac{1}{M} \sum\limits_{m=1}^{M}\nabla f\left(\hat{x}_{k,m}\right)\right\|^2 \\
    =& \frac{1}{M^2} \left\|\sum\limits_{m=1}^{M}\left(\nabla f(x_k)-
       \nabla f\left(\hat{x}_{k,m}\right)\right)\right\|^2\\
    \leq& \frac{1}{M} \sum\limits_{m=1}^{M}\left\|\left(\nabla f(x_k)- \nabla f\left(\hat{x}_{k,m}\right)\right)\right\|^2\\
    \leq & \frac{L_T^2}{M}\sum_{m=1}^M \|x_k - \hat{x}_{k,m}\|^2 \\
    \leq&L_T^2\max_{k\in \{1,\cdots, M\}}\|x_k - \hat{x}_{k,m}\|^2 \\
    =&L_T^2\|x_k - \hat{x}_{k,\mu}\|^2. \quad \text{(let $\mu:={\arg\max}_{m\in \{1,\cdots, M\}}\|x_k - x_{k-\tau_{k, m}}\|^2$)}
  \end{align*}
  It follows that
  \begin{align}
    T_1 \leq& L_T^2\left\|x_k - \hat{x}_{k, \mu}\right\|^2\nonumber\\
    =& L_T^2 \left\| \sum_{j\in J\left(k,\mu\right)} \left(x_{j+1} - x_j\right) \right\|^2\nonumber\\
    =& L_T^2 \gamma^2 \left\|\sum_{j\in J\left(k,\mu\right)}\sum\limits_{m=1}^{M}(G(\hat{x}_{j, m} ; \xi_{j,m}))_{i_j}e_{i_j}\right\|^2 \nonumber\\
    =& L_T^2 \gamma^2 \left\|\sum_{j\in J\left(k,\mu\right)}\sum\limits_{m=1}^{M} \left[(G(\hat{x}_{j, m} ; \xi_{j,m}))_{i_j}- \nabla_{i_j}f\left(\hat{x}_{j, m}\right) \right]e_{i_j} + \sum_{j\in J\left(k,\mu\right)} \sum\limits_{m=1}^{M}\nabla_{i_j}f\left(\hat{x}_{j, m}\right) e_{i_j}\right\|^2\nonumber \\
    \leq &2L_T^2 \gamma^2  \left( \underbrace{\left\|\sum_{j\in J\left(k,\mu\right)}\sum\limits_{m=1}^{M} \left[(G(\hat{x}_{j, m} ; \xi_{j,m}))_{i_j}- \nabla_{i_j}f\left(\hat{x}_{j, m}\right) \right]e_{i_j}\right\|^2}_{T_3} + \underbrace{\left\|\sum_{j\in J\left(k,\mu\right)} \sum\limits_{m=1}^{M}\nabla_{i_j}f\left(\hat{x}_{j, m}\right) e_{i_j}\right\|^2}_{T_4} \right)
           \label{eq:8}
  \end{align}
  where the last inequality uses the fact that
  $\|a + b\|^2 \le 2\|a\|^2 + 2\|b\|^2$ for any real vector $a$ and
  $b$. Taking the expectation in terms of $i_j$ and $\xi_{j,*}$ with
  all $j$'s in $J(k,\mu)$ for $T_3$, we have
  \begin{align}
    \nonumber &    \mathbb{E}_{\xi_{j,*}, i_j, j\in J(k, \mu)}(T_3) \\
    =&  \mathbb{E}_{\xi_{j,*}, i_j, j\in J(k, \mu)}\left( \left\|\sum_{j\in J(k,\mu)}\sum\limits_{m=1}^{M} \left[(G(\hat{x}_{j, m} ; \xi_{j,m}))_{i_j} - \nabla_{i_j}f(\hat{x}_{j, m})\right ]e_{i_j}\right\|^2\right)
       \nonumber \\
    =& \mathbb{E}_{\xi_{j,*}, i_j, j\in J(k, \mu)}\left( \sum_{j\in J(k,\mu)} \left\|\sum\limits_{m=1}^{M}\left[(G(\hat{x}_{j, m} ; \xi_{j,m}))_{i_j} - \nabla_{i_j}f(\hat{x}_{j, m})\right ]e_{i_j}\right\|^2 \right)
       \nonumber\\
              &+\mathbb{E}_{\xi_{j,*}, i_j, j\in J(k, \mu)}\Bigg(\sum_{j^{\prime\prime}\neq j^\prime, j^{\prime\prime}, j^\prime\in J(k,\mu)} \Bigg\langle \sum\limits_{m=1}^{M}\left[(G(\hat{x}_{j^{\prime\prime}, m} ; \xi_{j^{\prime\prime},m}))_{i_{j^{\prime\prime}}} - \nabla_{i_{j^{\prime\prime}}}f(\hat{x}_{j^{\prime\prime}, m})\right ]e_{i_{j^{\prime\prime}}},\nonumber\\
              &\sum\limits_{m=1}^{M}\left[(G(\hat{x}_{{j^\prime}, m} ; \xi_{j^\prime,m}))_{i_{j^\prime}} - \nabla_{i_{j^\prime}}f(\hat{x}_{{j^\prime}, m})\right ]e_{i_{j^\prime}}\Bigg\rangle \Bigg)
                \nonumber \\
    =&  \mathbb{E}_{\xi_{j,*}, i_j, j\in J(k, \mu)}\left( \sum_{j\in J(k,\mu)} \left\|\sum\limits_{m=1}^{M}\left[(G(\hat{x}_{j, m} ; \xi_{j,m}))_{i_j} - \nabla_{i_j}f(\hat{x}_{j, m})\right ]e_{i_j}\right\|^2 \right)
       \nonumber\\
    =&  {1\over n}\mathbb{E}_{\xi_{j,*},  j\in J(k, \mu)}\left( \sum_{j\in J(k,\mu)} \left\|\sum\limits_{m=1}^{M}\left[G(\hat{x}_{j, m} ; \xi_{j,m}) - \nabla f(\hat{x}_{j, m})\right ]\right\|^2 \right)
       \nonumber\\
    =& {1\over n}\mathbb{E}_{\xi_{j,*},  j\in J(k, \mu)}\left( \sum_{j\in J(k,\mu)} \sum\limits_{m=1}^{M}\left\|\left[G(\hat{x}_{j, m} ; \xi_{j,m}) - \nabla f(\hat{x}_{j, m})\right ]\right\|^2 \right)
       \nonumber\\
    \leq &
           \frac{ TM\sigma^2}{n}.
           \label{eq:7}
  \end{align}
  where the second last equality is due to \eqref{cross0} and the
  third equality is due to
  \begin{align*}
    &\mathbb{E}_{\xi_{{j}_*}, i_{j}, j \in J(k, \mu)}\Bigg(\sum_{j^{\prime\prime}\neq j^\prime, j^{\prime\prime}, j^\prime\in J(k,\mu)} \Bigg\langle \sum\limits_{m=1}^{M}\Bigg[(G(\hat{x}_{j^{\prime\prime}, m} ; \xi_{j^{\prime\prime},m}))_{i_{j^{\prime\prime}}} - \nabla_{i_{j^{\prime\prime}}}f(\hat{x}_{j^{\prime\prime}, m})\Bigg ]e_{i_{j^{\prime\prime}}},\\
    &\sum\limits_{m=1}^{M}\left[(G(\hat{x}_{{j^\prime}, m} ; \xi_{j^\prime,m}))_{i_{j^\prime}} - \nabla_{i_{j^\prime}}f(\hat{x}_{{j^\prime}, m})\right ]e_{i_{j^\prime}}\Bigg\rangle \Bigg)\\
    =  & 2\mathbb{E}_{\xi_{{j}_*}, i_{j}, j \in J(k, \mu)}\Bigg(\sum_{j^{\prime\prime}> j^\prime, j^{\prime\prime}, j^\prime\in J(k,\mu)} \Bigg\langle \sum\limits_{m=1}^{M}\Bigg[(G(\hat{x}_{j^{\prime\prime}, m} ; \xi_{j^{\prime\prime},m}))_{i_{j^{\prime\prime}}} - \nabla_{i_{j^{\prime\prime}}}f(\hat{x}_{j^{\prime\prime}, m})\Bigg ]e_{i_{j^{\prime\prime}}},\\
    &\sum\limits_{m=1}^{M}\left[(G(\hat{x}_{{j^\prime}, m} ;
      \xi_{j^\prime,m}))_{i_{j^\prime}} -
      \nabla_{i_{j^\prime}}f(\hat{x}_{{j^\prime}, m})\right
      ]e_{i_{j^\prime}}\Bigg\rangle \Bigg)\\
    = & 2\mathbb{E}_{\xi_{{j}_*}, i_{j}, j \in J(k, \mu)}\Bigg(\sum_{j^{\prime\prime}> j^\prime, j^{\prime\prime}, j^\prime\in J(k,\mu)} \Bigg\langle \sum\limits_{m=1}^{M}\Bigg[(\mathbb{E}_{\xi_{j^{\prime\prime}, m}}G(\hat{x}_{j^{\prime\prime}, m} ; \xi_{j^{\prime\prime},m}))_{i_{j^{\prime\prime}}} - \nabla_{i_{j^{\prime\prime}}}f(\hat{x}_{j^{\prime\prime}, m})\Bigg ]e_{i_{j^{\prime\prime}}},\\
    &\sum\limits_{m=1}^{M}\left[(G(\hat{x}_{{j^\prime}, m} ; \xi_{j^\prime,m}))_{i_{j^\prime}} - \nabla_{i_{j^\prime}}f(\hat{x}_{{j^\prime}, m})\right ]e_{i_{j^\prime}}\Bigg\rangle \Bigg)\\
    =& 0.
  \end{align*}

  Taking the expectation in terms of $i_j$ with all $j$'s in
  $J(k, \mu)$ for $T_4$, we have
  \begin{align}
    \nonumber
    &\mathbb{E}_{i_j, j\in J(k, \mu)}(T_4)  \\
    =& \mathbb{E}_{i_j, j\in J(k, \mu)}\left(  \left\|\sum_ {j\in J(k,\mu)} \sum\limits_ {m=1}^ {M}\nabla_ {i_j} f (\hat{x}_{j, m}) e_ {i_j}\right\|^2\right)
       \nonumber\\
    =& \mathbb{E}_{i_j, j\in J(k,\mu)}\Bigg[ \sum_{j\in J(k,\mu)}\left\| \sum\limits_{m=1}^{M}\nabla_{i_j}f(\hat{x}_{j, m}) e_{i_j}\right\|^2\nonumber\\
    &+ 2\sum_{j^{\prime\prime}> j^\prime, j^{\prime\prime}, j^\prime \in J(k,\mu)} \left\langle \sum\limits_{m=1}^{M}\nabla_{i_{j^{\prime\prime}}} f(\hat{x}_{j^{\prime\prime}, m})e_{i_{j^{\prime\prime}}}, \sum\limits_{m=1}^{M}\nabla_{i_{j^\prime}}f(\hat{x}_{j^\prime, m})e_{i_{j^\prime}} \right\rangle \Bigg]
      \nonumber\\
    =&  \mathbb{E}_{i_j, j\in J(k,\mu)}\Bigg[ {1\over n}\sum_{j\in J(k,\mu)}\left\| \sum\limits_{m=1}^{M}\nabla f(\hat{x}_{j, m}) \right\|^2  \nonumber\\
    & +  2 \sum_{j^{\prime\prime}> j^\prime, j^{\prime\prime}, j^\prime \in J(k,\mu)}  \left\langle \sum\limits_{m=1}^{M}\nabla_{i_{j^{\prime\prime}}} f(\hat{x}_{j^{\prime\prime}, m})e_{i_{j^{\prime\prime}}}, \sum\limits_{m=1}^{M}\nabla_{i_{j^\prime}}f(\hat{x}_{j^\prime, m})e_{i_{j^\prime}} \right\rangle \Bigg]
      \nonumber \\
    \leq &
           \mathbb{E}_{i_j, j\in J(k,\mu)}\Bigg[ {1\over n}\sum_{j\in J(k,\mu)}\left\| \sum\limits_{m=1}^{M}\nabla f(\hat{x}_{j, m}) \right\|^2   \nonumber\\
    &+  {1\over n} \sum_{j^{\prime\prime} > j^\prime, j^{\prime\prime}, j^\prime \in J(k,\mu)} \left( {1\over \alpha}\left\|\sum\limits_{m=1}^{M}\nabla f(\hat{x}_{j^{\prime\prime}, m})\right\|^2 + {\alpha \over n} \left\|\sum\limits_{m=1}^{M}\nabla f(\hat{x}_{j^\prime, m})\right\|^2\right)\Bigg] \text{  (Let }\alpha=\sqrt{n}\text{)}
      \nonumber\\
    \leq &
           \left({\sqrt{n}+T-1\over n^{3/2}}\right)\sum_{j\in
           J(k,\mu)}  \mathbb{E}_{i_j, j\in
           J(k,\mu)}\left(\left\|\sum\limits_{m=1}^{M}\nabla
           f(\hat{x}_{j, m})\right\|^2\right),
           \label{eq:6}
  \end{align}
  where the second last inequality is due to that for any $j^{\prime\prime}> j^\prime$:
  \begin{align*}
    &\mathbb{E}_{i_{j^{\prime\prime}}, i_{j^\prime}}\left\langle \sum\limits_{m=1}^{M}\nabla_{i_{j^{\prime\prime}}} f\left(\hat{x}_{j^{\prime\prime}, m}\right)e_{i_{j^{\prime\prime}}}, \sum\limits_{m=1}^{M}\nabla_{i_{j^\prime}}f\left(\hat{x}_{j^\prime, m}\right)e_{i_{j^\prime}} \right\rangle
    \\
    =&  {1\over n}\mathbb{E}_{i_{j^\prime}} \left\langle \sum\limits_{m=1}^{M}\nabla f\left(\hat{x}_{j^{\prime\prime}, m}\right), \sum\limits_{m=1}^{M}\nabla_{i_{j^\prime}}f\left(\hat{x}_{j^\prime, m}\right)e_{i_{j^\prime}} \right\rangle
    \\
    \leq &
           {1\over n} \mathbb{E}_{i_{j^\prime}}\left( {1\over 2\alpha}\left\|\sum\limits_{m=1}^{M}\nabla f\left(\hat{x}_{j^{\prime\prime}, m}\right)\right\|^2 + {\alpha \over 2} \left\|\sum\limits_{m=1}^{M}\nabla_{i_{j^\prime}}f\left(\hat{x}_{j^\prime, m}\right)e_{i_{j^\prime}}\right\|^2\right)
    \\
    =&
       {1\over n} \mathbb{E}_{i_{j^\prime}}\left( {1\over 2\alpha}\left\|\sum\limits_{m=1}^{M}\nabla f\left(\hat{x}_{j^{\prime\prime}, m}\right)\right\|^2 + {\alpha \over 2n} \left\|\sum\limits_{m=1}^{M}\nabla f\left(\hat{x}_{j^\prime, m}\right)\right\|^2\right).
  \end{align*}
  Take full expectation on both sides of \eqref{eq:6}, \eqref{eq:7}
  and \eqref{eq:8}. Then substituting the upper bound of
  $\mathbb{E}(T_3)$ and $\mathbb{E}(T_4)$ into $\mathbb{E}(T_1)$, we have
  \begin{align}
    \label{eq:9}
    \mathbb{E}(T_1) \le 2L_T^2\gamma^2 \left( \frac{TM\sigma^2}{n} +\left({\sqrt{n}+T-1\over n^{3/2}}\right)\sum_{j\in
    J(k,\mu)}  \mathbb{E}\left(\left\|\sum\limits_{m=1}^{M}\nabla
    f(\hat{x}_{j, m})\right\|^2\right) \right)
  \end{align}

  Take full expectation on both sides of \eqref{eq:10} and
  \eqref{eq:lip}. Substituting $\mathbb{E}(T_1)$ and $\mathbb{E}(T_2)$
  into \eqref{eq:lip}, we have
  \begin{align}
    \mathbb{E} \left(f\left(x_{k+1}\right)\right) \leq& \mathbb{E}\left(   f \left(x_k\right)\right) - \frac{M\gamma}{2n}  \Bigg[\mathbb{E}\left(  \left\|\nabla f \left(x_k\right)\right\|^2\right) + \mathbb{E} \left(  \left\|\frac{1} {M}\sum\limits_ {m=1}^ {M}\nabla f \left(\hat{x}_{k, m}\right)\right\|^2\right) \nonumber\\
                                                      &- 2L_T^2 \gamma^2 \left(\frac{ TM\sigma^2}{n} + \left({\sqrt{n}+T-1\over n^{3/2}}\right)\sum_{j\in J\left(k,\mu\right)}  \mathbb{E}\left(\left\|\sum\limits_{m=1}^{M}\nabla f\left(\hat{x}_{j, m}\right)\right\|^2\right)\right)  \Bigg]\nonumber\\
                                                      &+\frac{\Lmax\gamma^2}{2n}\left(M \sigma^2 + \mathbb{E}\left(   \left\|\sum_ {m=1}^ M \nabla f \left(\hat{x}_{k, m}\right) \right\|^2\right)\right)\nonumber\\
    =&  \mathbb{E}\left(   f \left(x_k\right)\right) - \frac{M\gamma}{2n} \mathbb{E}\left(  \left\|\nabla f \left(x_k\right)\right\|^2\right)  - \frac{\gamma}{2Mn} \mathbb{E}\left(  \left\|\sum\limits_ {m=1}^ {M}\nabla f \left(\hat{x}_{k, m}\right)\right\|^2\right)\nonumber\\
                                                      &+ \frac{M\gamma^3}{n} L_T^2 \left(\left({\sqrt{n}+T-1\over n^{3/2}}\right)\sum_{j\in J\left(k,\mu\right)}  \mathbb{E}\left(\left\|\sum\limits_{m=1}^{M}\nabla f\left(\hat{x}_{j, m}\right)\right\|^2\right)\right) \nonumber\\
                                                      &+\frac{\Lmax\gamma^2}{2n}\left(\mathbb{E} \left\|\sum_{m=1}^M \nabla f\left(\hat{x}_{k, m}\right) \right\|^2\right) + \frac{L_T^2TM^2\gamma^3}{n^2} \sigma^2 + \frac{\Lmax M\gamma^2}{2n} \sigma^2.
                                                        \label{eq:k1k}
  \end{align}
  Summarizing \eqref{eq:k1k} from $k = 1$ to $K$, we have
  \begin{align*}
    \mathbb{E} \left(f\left(x_{K+1}\right)\right) \leq& f\left(x_1\right) - \sum_{k=1}^K\frac{M\gamma}{2n} \mathbb{E}\left(  \left\|\nabla f \left(x_k\right)\right\|^2\right)  - \frac{\gamma}{2Mn} \sum_{k=1}^K\mathbb{E}\left(  \left\|\sum\limits_ {m=1}^ {M}\nabla f \left(\hat{x}_{k, m}\right)\right\|^2\right)\\
                                                      &+ \frac{M\gamma^3}{n} L_T^2 \left(\left({\sqrt{n}+T-1\over n^{3/2}}\right)\sum_{k=1}^K\sum_{j\in J\left(k,\mu\right)}  \mathbb{E}\left(\left\|\sum\limits_{m=1}^{M}\nabla f\left(\hat{x}_{j, m}\right)\right\|^2\right)\right) \\
                                                      &+\frac{\Lmax\gamma^2}{2n}\sum_{k=1}^K\left(\mathbb{E}\left(   \left\|\sum_ {m=1}^ M \nabla f \left(\hat{x}_{k, m}\right) \right\|^2\right)\right) + \frac{KL_T^2TM^2\gamma^3}{n^2} \sigma^2 + \frac{K\Lmax M\gamma^2}{2n} \sigma^2 \\
    \le &f\left(x_1\right) - \sum_{k=1}^K\frac{M\gamma}{2n} \mathbb{E}\left(  \left\|\nabla f \left(x_k\right)\right\|^2\right)  - \frac{\gamma}{2Mn} \sum_{k=1}^K\mathbb{E}\left(  \left\|\sum\limits_ {m=1}^ {M}\nabla f \left(\hat{x}_{k, m}\right)\right\|^2\right)\\
                                                      &+ \frac{MT\gamma^3}{n} L_T^2 \left(\left({\sqrt{n}+T-1\over n^{3/2}}\right)\sum_{k=1}^K\mathbb{E}\left(\left\|\sum\limits_{m=1}^{M}\nabla f\left(\hat{x}_{k, m}\right)\right\|^2\right)\right) \\
                                                      &+\frac{\Lmax\gamma^2}{2n}\sum_{k=1}^K\left(\mathbb{E}\left(   \left\|\sum_ {m=1}^ M \nabla f \left(\hat{x}_{k, m}\right) \right\|^2\right)\right) + \frac{KL_T^2TM^2\gamma^3}{n^2} \sigma^2 + \frac{K\Lmax M\gamma^2}{2n} \sigma^2 \\
    \le &f\left(x_1\right) - \sum_{k=1}^K\frac{M\gamma}{2n} \mathbb{E}\left(  \left\|\nabla f \left(x_k\right)\right\|^2\right) + \frac{KL_T^2TM^2\gamma^3} {n^2} \sigma^ 2 + \frac{K\Lmax M\gamma^2}{2n} \sigma^2 \\
                                                      &+
                                                        \left( \frac{MT L_T^2\left(\sqrt{n}+T-1\right)\gamma^3}{n^{5/2}} - \frac{\gamma}{2Mn} +\frac{\Lmax\gamma^2}{2n}\right)\sum_{k=1}^K\left(\mathbb{E} \left\|\sum_{m=1}^M \nabla f\left(\hat{x}_{k, m}\right) \right\|^2\right) \\
    \le & f\left(x_1\right) - \sum_{k=1}^K\frac{M\gamma}{2n} \mathbb{E}\left(  \left\|\nabla f \left(x_k\right)\right\|^2\right) + \frac{KL_T^2TM^2\gamma^3}{n^2} \sigma^2 + \frac{K\Lmax M\gamma^2}{2n} \sigma^2
  \end{align*}
  where the last inequality comes from \eqref{alg2ass}. Together with
  $\mathbb{E}(f(x_{k+1}))\ge f(x^*)$, we have
  \begin{align}
    \frac{1}{K}\sum_{t=1}^K \mathbb{E}\|\nabla f(x_t)\|^2\leq  \frac{2n}{KM\gamma} (f(x_1) - f(x^*))+ \frac{2L_T^2TM\gamma^2}{n} \sigma^2 + \Lmax\gamma \sigma^2.
  \end{align}
  It completes the proof.
\end{proof}

{\noindent \bf Proofs to Corollary~\ref{coro:inconread}}
\begin{proof}
  From the definition of the steplength \eqref{eq:coro_11} and the
  lower bound of $K$ \eqref{eq:coro_22}, we have
  \begin{align}
    \nonumber
    \gamma &= \sqrt{\frac{n(f(x_1) - f(x^*))}{L_TM\sigma^2}}\frac{1}{\sqrt{K}}\\
    \nonumber
           &\le  \sqrt{\frac{n(f(x_1) - f(x^*))}{L_TM\sigma^2}}\sqrt{\frac{\sqrt{n} \sigma ^2}{16 (f(x_1) - f(x^*)) L_T M \left(n^{3/2}+4 T^2\right)}}\\
           &= \frac{1}{4}\sqrt{\frac{1}{n^{3/2}+4 T^2}}\frac{n^{3/4}}{L_TM}
             \label{eq:proof:alg2:1}
  \end{align}
  which gives an upper bound for $\gamma$. We can further relax this
  upper bound by
  \begin{align}
    \gamma \leq \frac{1}{4}\sqrt{\frac{1}{n^{3/2}+4 T^2}}\frac{n^{3/4}}{L_TM} \leq \frac{1}{4}\sqrt{\frac{1}{4 T^2}}\frac{n^{3/4}}{L_TM} < \frac{2n}{L_TMT}.
    \label{eq:proof:alg2:2}
  \end{align}
  Next we will show that the steplength satisfies the condition in
  \eqref{alg2ass}:
  \begin{align*}
    \text{LHS}= &\frac{2M^2T L_T^2\left(\sqrt{n}+T-1\right)\gamma^2}{n^{3/2}} +2M\Lmax\gamma \\
    =& \frac{2M^2T L_T^2\left(\sqrt{n}+T-1\right)}{n^{3/2}}\frac{1}{16}\frac{1}{n^{3/2}+4 T^2}\frac{n^{3/2}}{L_T^2M^2}
       +\frac{\Lmax M}{2} \sqrt{\frac{1}{n^{3/2}+4 T^2}}\frac{n^{3/4}}{L_TM}\\
    \le& \frac{1}{8}\frac{T \left(\sqrt{n}+T\right)}{n^{3/2}+4 T^2} + \frac{1}{2}\\
    \le& \frac{1}{8}\frac{T\sqrt{n} + T^2}{n^{3/2}+4 T^2} + \frac{1}{2}\\
    \le &\frac{1}{8}\frac{3T^2}{n^{3/2}+4 T^2} + \frac{1}{8}\frac{2n}{n^{3/2}+4 T^2} + \frac{1}{2}\\
    \le& \frac{1}{8}\cdot\frac{3}{4} + \frac{1}{8}\cdot\frac{1}{2} + \frac{1}{2}\\
    \le & 1 = \text{RHS},
  \end{align*}
  where the first inequality uses the fact $\Lmax \leq L_T$ and the
  first inequality uses the upper bound of $\gamma$ in
  \eqref{eq:proof:alg2:1}. The third last
  inequality comes from $\sqrt{n}T\le2n+2T^2$. It means that the
  condition \eqref{alg2ass} in Theorem~\ref{thm:algo2} is satisfied
  globally. Then we can safely apply \eqref{eq:algo2} in
  Theorem~\ref{thm:algo2}:
  \begin{align*}
    \frac{1}{K}\sum_{t=1}^K \mathbb{E}\left\|\nabla f\left(x_t\right)\right\|^2\leq&  \frac{2n}{KM\gamma} \left(f\left(x_1\right) - f\left(x^*\right)\right)+ \frac{2L_T^2TM\gamma^2}{n} \sigma^2 + \Lmax\gamma \sigma^2\\
    \le &\frac{2n}{KM\gamma} \left(f\left(x_1\right) - f\left(x^*\right)\right)+ 5 L_T \gamma \sigma^2 \\
    =& \frac{6 \sqrt{2 \left(f\left(x_1\right) - f\left(x^*\right)\right) L_T n } \sigma}{\sqrt{KM}}
  \end{align*}
  where the second inequality is due to the upper bound of $\gamma$ in
  \eqref{eq:proof:alg2:2} and the last equality is acquired by
  substituting $\gamma$ by its definition in \eqref{eq:coro_11}. It
  completes the proof.
\end{proof}

\end{document}